\documentclass[a4paper,12pt]{article}

\usepackage{latexsym}
\usepackage{amssymb}
\usepackage{amsthm}
\usepackage{amsmath}
\usepackage{amssymb}
\usepackage{enumerate}
\usepackage{color}
\usepackage{dsfont} 

\newcommand{\p}{\partial}
\newcommand{\eps}{\varepsilon}

\newcommand{\dd}{\hspace{1pt}{\rm d}\hspace{0.0pt}}
\newcommand{\ee}{{\rm e}\hspace{1pt}}

\newcommand{\R}{\mathbb R}
\newcommand{\N}{\mathbb N}
\newcommand{\C}{\mathbb C}

\newcommand{\cB}{{\mathcal B}}
\newcommand{\cC}{{\mathcal C}}
\newcommand{\cD}{{\mathcal D}}
\newcommand{\cE}{{\mathcal E}}

\newcommand{\cG}{{\mathcal G}}

\newcommand{\cM}{{\mathcal M}}
\newcommand{\cN}{{\mathcal N}}
\newcommand{\cO}{{\mathcal O}}

\newcommand{\cS}{{\mathcal S}}

\newcommand{\rE}{\ensuremath{\mathrm{E}}}
\newcommand{\rV}{\ensuremath{\mathrm{V}}}

\newcommand{\Cinf}{\ensuremath{{\mathcal C}^\infty}}

\newcommand{\EM}{\ensuremath{{\cE}_{\mathrm{M}}}}

\newcommand{\supp}{\mathop{\mathrm{supp}}}
\newcommand{\esssup}{\mathop{\mathrm{ess\,sup}}}

\newcommand{\eins}{\mathds{1}}

\theoremstyle{plain}
\newtheorem{thm}{Theorem}[section]
\newtheorem{lemma}{Lemma}[section]
\newtheorem{prop}{Proposition}[section]
\newtheorem{cor}{Corollary}[section]
\theoremstyle{definition}
\newtheorem{ex}{Example}[section]
\newtheorem{rem}{Remark}[section]
\newtheorem{defn}{Definition}[section]

\begin{document}

\title{Generalized solutions to hyperbolic systems with random field coefficients}

\vspace{5mm}
\author{Jelena Karaka\v{s}evi\'{c}\footnote{Electronic mail: jelena.karakasevic.mail@gmail.com}, Michael Oberguggenberger\footnote{Electronic mail: michael.oberguggenberger@uibk.ac.at},
        Martin Schwarz\footnote{Electronic mail: math@mschwarz.eu}\\
        {\small Unit for Engineering Mathematics}\\
         {\small University of Innsbruck, Technikerstr.\ 13, A-6020 Innsbruck, Austria}
         }
\date{}         
\maketitle

{\bf Abstract:}
The paper addresses linear hyperbolic systems in one space dimension with random field coefficients. In many applications, a low degree of regularity of the paths of the coefficients is required, which is not covered by classical stochastic analysis. For this reason, we place our analysis in the framework of Colombeau algebras of generalized functions. We obtain new characterizations of Colombeau stochastic processes and establish existence and uniqueness of
solutions in this framework. A number of applications to stochastic wave and transport equations are given and the Colombeau solutions are related to classical weak solutions, when the latter exist.
\vspace{5mm}
\vspace{5mm}
%
%
\section{Introduction}
\label{Sec:intro}
%
%
This paper is devoted to the Cauchy problem for linear hyperbolic systems of the form of
	\begin{equation}
	\begin{array}{rcl}
	\left(\partial_t+\lambda(x,t)\partial_x\right)u&=&f(x,t)u+g(x,t), \quad(x,t)\in\R^2,\\
	u(x,0)&=&u_0(x), \quad x\in\R,
	\end{array}
	\label{eq:hsys}
	\end{equation}
where $u=(u_1,\ldots,u_n)$, $g=(g_1,\ldots,g_n)$; $\lambda$ and $f$ are $(n\times n)$-matrix valued functions. The coefficient matrix $\lambda$ is assumed to be real-valued and diagonal,
$\lambda = {\rm diag}\,(\lambda_1.\ldots,\lambda_n)$.
The focus of the paper will be on the case where the coefficients as well as the initial data are random fields on some probability space $(\Omega,\Sigma,P)$.

Classical random solutions can be found, provided the paths of the random fields are sufficiently regular. For example, if the paths of $\lambda$ are continuous in $(x,t)$ and Lipschitz with respect to the first variable, while the entries
of $f$ and $g$ are continuous, the corresponding system of integral equations has an almost surely unique pathwise solution $u$. This assertion already fails in the deterministic case if the paths of $\lambda$ are merely H\"{o}lder continuous
\cite[Example 6.1]{Flandoli:2010}.

From the viewpoint of applications it is desirable to admit coefficients of even lower regularity. The main applicative interest lies in wave propagation and transport in randomly perturbed media. In engineering applications random material parameters are often modelled by random fields with continuous, but nowhere differentiable paths \cite{Ghanem:1991}. A frequently used example is the Ornstein-Uhlenbeck process \cite[Section 8.3]{Arnold:1974} whose autocorrelation function is continuous, but not differentiable \cite[Formula (2.39)]{Ghanem:1991}. Another source of examples arises from randomly layered media \cite{Fouque:2007}. Certain scaling limits lead to diffusion processes as coefficients, such as Brownian motion
\cite[Section 5.2]{Fouque:2007}. Scaling limits in the so-called Goupillaud medium may give rise to arbitrary L\'evy processes \cite{BOS2017}. It may even be necessary to insert distributional derivatives of such processes in the coefficient functions, see Example\;\ref{ex:wavegeometric} below, when the geometric wave equation in a randomly perturbed geometry is studied.

In short, we wish to admit generalized stochastic processes as coefficients (or initial data) in \eqref{eq:hsys}. This entails that a meaning has to be given to the products of the coefficients with the generalized solution in \eqref{eq:hsys}.
To this end, we adopt the framework of Colombeau algebras of generalized functions. These are differential algebras of generalized functions which contain the space of Schwartz distributions as a subspace. They have found a wealth of applications in nonlinear partial differential equations and in nonsmooth geometry \cite{GKOS}. Colombeau stochastic processes have been considered since the 1990s \cite{AHR,MO1995}. We will rely on the notion of Colombeau stochastic processes as developed in \cite{GOPS2018,Mirkov:2009}, with minor modifications.

Colombeau stochastic processes or Colombeau random fields are defined as equivalence classes of families $(u_\eps)_{\eps\in(0,1]}$ of random processes with smooth paths which satisfy certain asymptotic properties as $\eps\to 0$. At fixed $\eps$, \eqref{eq:hsys} turns
into a system of random differential equations, which can be solved by classical techniques \cite{Bunke:1972,Kloeden:2017}. The key ingredient of the existence and uniqueness theory is to prove the asymptotic estimates for the representatives
$(u_\eps)_{\eps\in(0,1]}$.

Section\;\ref{Sec:classdet} briefly summarizes the classical deterministic theory for constructing continuous or differentiable solutions to \eqref{eq:hsys} by integrating along characteristic curves. We also recall a basic a priori estimate which plays a central role in deriving the asymptotic bounds required in later sections. In Section\;\ref{Sec:classrand} random classical solutions to \eqref{eq:hsys} are constructed, and the required conditions on the coefficients are highlighted.
Though the results of these two sections are well-known, we found it useful to briefly recall them because the estimates become important on a more general level in the construction of the Colombeau solutions.

Section\;\ref{Sec:Colombeaurandfun} starts with recalls on deterministic Colombeau theory and then presents spaces of Colombeau stochastic process (with or without moments). It is indicated how random Schwartz distributions can be embedded. A central result is a new asymptotic characterization (Lemma\;\ref{lemma:supequiv}) which forms a decisive step beyond \cite{GOPS2018,Mirkov:2009}. The existence of expectation values, moments and the autocovariance function as members of Colombeau algebras of generalized functions is discussed.

The main new results are in Section\;\ref{Sec:exis}: existence and uniqueness of Colombeau solutions to the stochastic system \eqref{eq:hsys} without moments (Theorem\;\ref{thm:gensolu0}) and with moments (Theorem\;\ref{thm:gensolup}).
It is also shown how the generalized solutions relate to the classical solutions from Section\;\ref{Sec:classrand}, when the latter exist. As a special case, generalized solutions to the scalar transport equation are considered. It is shown that
solutions exist under weaker assumptions, provided the drift coefficient depends on one of the variables $(x,t)$ only.

Section\;\ref{Sec:apps} is devoted to applications. We start out with the one-dimensional wave equation with lower order terms (Example\;\ref{ex:wavegeneral}) and its reduction to a first order system. This is applied to the stochastic wave equation with random propagation speed depending on $x$ only (Example\;\ref{ex:wave}). We turn to the constant coefficient wave equation with additive noise in the Colombeau sense in Example\;\ref{ex:wavenoise}. If the noise is embedded space-time Gaussian white noise, we show that the unique Colombeau solution is associated (in the limit as $\eps\to 0$) to the well-known classical random field solution \cite{Walsh:1984}. Example\;\ref{ex:wavegeometric} addresses the geometric wave equation on the graph of a smooth function, viewed as a one-dimensional Riemannian manifold. The geometric wave equation fully decouples in a pair of transport equations. If the graph is perturbed by a stochastic process with continuous paths of locally infinite length, there appears to be no applicable classical solution concept, but Colombeau solutions can still be constructed. It is shown that in this case waves do not propagate at all, while for continuously differentiable perturbations, the propagation is as usual.
The final Example\;\ref{ex:Ogawa} takes up Ogawa's analysis \cite{Ogawa} of of the transport equation with white noise (in time) as coefficient. We construct a unique Colombeau solution and recover Ogawa's results in the limit as $\eps\to 0$.
Note that the effect of adding time-dependent noises to the drift has become an issue of considerable interest in the past decade \cite{Flandoli:2019,Flandoli:2011,Flandoli:2010}.

The Appendix collects some facts about embedding distributions into Colom\-beau algebras. Working on $\R^d$, an especially simple form of embedding can be used. The form of this embedding was introduced in \cite[Section 12]{MOBook} without much detail. For the purpose of applying it to generalized stochastic processes, we provide the required estimates and also show that it coincides with the usual embedding as worked out in \cite[Section 1.2.2]{GKOS}. Following recent clarifications \cite{MO:2024}, generalized stochastic processes -- random Schwartz distributions or $L^p$-valued random linear functionals -- can be embedded into spaces of Colombeau stochastic processes in a simple and unified way.
Examples of such processes as used in Section\;\ref{Sec:apps} are given, and it is shown how the requirements of the existence theory from Section\;\ref{Sec:exis} can be met. The Appendix ends with two counterexamples showing that certain relations between different variants of Colombeau spaces addressed in Section\;\ref{Sec:Colombeaurandfun} do not hold.

We conclude the introduction with a short survey of the existing literature. In the setting of classical stochastic analysis there are few papers addressing hyperbolic equations with random or stochastic coefficients. The main focus has been on linear and nonlinear partial differential equations driven by noise terms of various degrees of non-regularity. For wave equations, standard references are \cite{DaPratoZabczyk:2014,DalangBook:2009,Walsh:1984}. Here is a selection of further references to hyperbolic systems (with deterministic coefficients) perturbed by noise: \cite{Coriasco:2019,Coriasco:2020,CarmonaNualart:1988,Kim:2000,Kim:2011,Marcus1991,Orsingher:1982,PeszatZabczyk:2000}.

Hyperbolic systems with random field coefficients have been considered in \cite{ChowHuang,Sukys:2016,Sukys:2012}
and for hyperbolic conservation laws with random flux \cite{Mishra:2016}. In all these papers, a certain degree of regularity is required for the coefficients.
Hyperbolic systems with highly singular coefficients have been considered in \cite{Coriasco:2021} by means of chaos expansion and the Wick product.
The case of scalar transport equations deserves special mention -- much more singular drift coefficients have been treated, especially containing time-dependent noise.
Starting with Ogawa's article \cite{Ogawa} on a transport equation with white noise in time added to the drift, a number of authors have considered more general transport equations with drift perturbed in time: \cite{Flandoli:2019,Flandoli:2011,Flandoli:2010,Neves:2022,Olivera:2019}, some of which were already mentioned above. Quasilinear first order equations were addressed in \cite{Kunita:1984}.
These type of equations with highly singular drift in time were intensively studied using rough paths: \cite{Catellier:2016,Friz:2016,Nilssen:2020}.
These concepts have also led to new results for deterministic transport equations with irregular coefficients \cite{BailleulGubinelli:2017,GaleatiGubinelli:2022}.

We end with relevant literature in Colombeau theory. In the deterministic case, existence and
uniqueness of solutions in Colombeau algebras for systems of type (\ref{eq:hsys}) with coefficients in Colombeau algebras was established in \cite{MO1989}, for symmetric hyperbolic systems in higher space dimensions in \cite{LafonMO}, for scalar hyperbolic pseudodifferential systems with Colombeau symbols in \cite{Hoer2004} (and with less restrictive conditions on the coefficients in \cite{DO2016}), for strictly hyperbolic pseudodifferential systems and higher order equations with Colombeau symbols in \cite{GarettoMO2015}. In the stochastic case most papers addressed linear and nonlinear SPDEs perturbed by a noise term given as a Colombeau stochastic process, but with deterministic coefficients. The papers
\cite{AHR,NedeljkovRajter:2002b,NedeljkovRajter:2002a,MORus1998b} addressed Colombeau solutions to nonlinear wave equations perturbed by a Colombeau stochastic process. So far, stochastic coefficients have been considered only for the Dirichlet problem (for an elliptic equation) in \cite{Pilipovic:2010}. In this paper, based on \cite{Selesi:2008}, the authors develop a theory of Colombeau stochastic processes with values in Kondratiev spaces, that is, a Colombeau variant of white noise analysis.
%
%
\section{Classical deterministic theory}
\label{Sec:classdet}
%
%
In this section we briefly recall the classical deterministic theory. The characteristic curves of system \eqref{eq:hsys} are the integral curves of the vector fields $\p_t + \lambda_i(x,t)\p_x$, where $\lambda_i$ is the $i$-th entry in the diagonal matris $\lambda$.
We adjust them so that they pass through a given point $x_0$ at time
$t = t_0$. Hence the characteristic curves are the solutions to
\begin{equation}
\begin{array}{rcl}
	\frac{\dd}{\dd t}\gamma_i(x_0,t_0,t)&=&\lambda_i\big(\gamma_i(x_0,t_0,t),t\big),\\
	\gamma_i(x_0,t_0,t_0)&=&x_0.
\end{array}
\label{eq:char}
\end{equation}
There are various sufficient conditions under which the characteristic curves exist globally and depend continuously on the starting values $(x_0,t_0)$. Anticipating what is needed in the random case, we will use the following condition:
\begin{itemize}
\item[(LB)] (Local Lipschitz condition and global bound) For every compact subset $K\subset \R$ and every compact time interval $I$, there are constants $\ell > 0$, $c>0$ such that, for all $i = 1,\ldots, n$,
\[
   |\lambda_i(x,t) - \lambda_i(y,t)| \leq \ell|x-y|\quad\mbox{for\ all\ } x,y\in K, t\in I,
\]
and
\[
   |\lambda_i(x,t)| \leq c \quad\mbox{for\ all\ } x\in \R, t\in I.
\]
\end{itemize}
\begin{lemma}\label{lemma:char}
Assume that the $\lambda_i = \lambda_i(x,t)$ are continuous functions of $(x,t)$ and that condition (LB) holds. Then equation \eqref{eq:char} has a unique solution which exists globally in time and depends continuously on the initial data $(x_0,t_0)$.

If in addition the coefficients $\lambda_i$ belongs to $\cC^k(\R^2)$, $k\geq 1$, the same is true of the characteristic flow, that is, $\gamma_i$ belongs to $\cC^k(\R^3)$.
\end{lemma}
\emph{Proof.} Rephrase problem \eqref{eq:char} as the integral equation
\begin{equation}
 \gamma_i(x_0,t_0,t) = x_0 - \int_t^{t_0}\lambda_i\big(\gamma_i(x_0,t_0,s),s\big)\dd s
 \label{eq:charint}
\end{equation}
and consider $(x_0,t_0)$ as fixed first, contained in the interior of some compact set $K\times I$. Let
\[
\cB_\tau = \{f\in\cC[t_0-\tau,t_0+\tau]:\sup_{t\in[t_0-\tau,t_0+\tau]}|f(t) - x_0| \leq 1\}.
 \]
The integral operator
\[
   (\cM h)(t) = x_0 -\int_t^{t_0} \lambda_i(h(s),s)\dd s
\]
defines a contraction on $\cB_\tau$ for sufficiently small $\tau$, thus producing a local solution $h(t) = \gamma(x_0,t_0,t)$ in $[t_0-\tau,t_0+\tau]$. The time of existence $\tau$ depends only on the constants $\ell$, $c$ in condition (LB) (and the choice of the compact sets $K$, $I$). Hence the solution can be continued globally in time.

Continuous dependence on the initial data $(x_0,t_0)$ can be inferred by Gronwall's inequality. The differentiability with respect to $(x_0,t_0,t)$ in case the coefficients $\lambda_i$ belong to $\cC^k(\R^2)$ follows, e.g., from
\cite[Corollary 4.1]{Hartman:1964}.
\qed

Returning to system \eqref{eq:hsys} and integrating along the characteristic curves results in the system of integral equations
\begin{align}
	u_i(x,t) & = u_{0i}\big(\gamma_i(x,t,0)\big)\label{eq:intsys}\\
                  &+ \int_0^t\Big(\sum_{j=1}^n f_{ij}\big(\gamma_i(x,t,\tau),\tau\big)u_j\big(\gamma_i(x,t,\tau),\tau\big) + g_i\big(\gamma_i(x,t,\tau),\tau)\big)\Big)\dd\tau \nonumber
\end{align}
for $i=1,\ldots,n$. We briefly sketch the existence theory for the system of integral equations \eqref{eq:intsys}.

Let $K_0 = [-\kappa,\kappa]$ be a compact interval. For $-T\leq t \leq T$, the region $K_T$ is defined by
\begin{equation}
  K_T = \{(x,t)\in \R^2 : -T\leq t\leq T, |x| \leq \kappa - c|t| \}
\label{eq:KT}
\end{equation}
where $c$ is the constant in (LB) associated with the interval $I = [-T,T]$.
The region $K_T$ is a domain of determinacy; for any point $(x_0,t_0)\in K_T$, all $n$ characteristic curves joining this point with the $x$-axis remain in $K_T$.
\begin{prop}\label{prop:classsolu}
Let the coefficients in $\lambda$, $f$ and $g$ belong to $\cC(\R^2)$. Assume further that $\lambda$ satisfies condition (LB). Then the following holds:

(a) If $u_0\in\cC(K_0)$, there is $T>0$ such that the problem (\ref{eq:intsys}) has a solution $u\in\cC(K_T)$.

(b) For whatever $T_0>0$, there is at most one solution $u\in\cC(K_{T_0})$.

(c) If $u_0\in\cC(\R)$, then problem (\ref{eq:intsys}) has a unique global solution $u\in\cC(\R^2)$.

(d) If $\lambda, f, g$ belong to $\cC^k(\R^2)$ and $u_0\in\cC^k(\R)$ for some $k \geq 1$, then $u\in\cC^k(K_{T_0})$.
\end{prop}
\emph{Sketch of proof.} Assertion (a) is obtained by showing that the integral operator on the right-hand side of \eqref{eq:intsys} defines a contraction in the ball of radius one in $\cC(K_T)$ around the free solution
\begin{equation}\label{eq:data}
  \Big(u_{0i}\big(\gamma_i(x,t,0)\big)+ \int_0^t g_i\big(\gamma_i(x,t,\tau),\tau)\big)\Big)\dd\tau\Big)_{i = 1,\ldots, n}
\end{equation}
for sufficiently small $T$. It is important to point out that the time of local existence depends only on the supremum norm of the $f_{ij}$, $g_i$ on $K_T$ and $u_0$ on $K_0$.

Point (b) is obtained by applying Lemma\;\ref{lemma:basiclemma} below to the difference of two solutions with the same data (which satisfies (\ref{eq:intsys}) with $g\equiv 0$, $u_0\equiv 0$).

(c) The global solvability follows from the a priori estimate implied by Lemma \ref{lemma:basiclemma} and the fact that $K_0$ can be chosen arbitrarily large.

(d) Differentiability is obtained by first applying the fixed point theorem in $\cC^k(K_T)$ for small $T$ and deriving a priori estimates on the derivatives using the lemma inductively. \qed

Finally, observe that for solutions in $\cC^1(K_T)$, the system of integral equations (\ref{eq:intsys}) is equivalent with the system of differential equations \eqref{eq:hsys}. This equivalence holds more generally for locally integrable solutions in the weak sense. Further details on the classical theory can be found in \cite[Section 13]{MOBook}.

In the sequel, we denote the maximum norm on vectors $g$ and the corresponding matrix norm on matrices $f$ simply by
\begin{equation}\label{eq:norms}
|g| = \max_{1\leq j \leq n}|g_j|, \quad |f| = \max_{1\leq i \leq n}\sum_{j=1}^n|f_{ij}|.
\end{equation}
\begin{lemma}\label{lemma:basiclemma}
Under the hypothesis of Proposition\; \ref{prop:classsolu}, assume that $u\in\cC(K_T)$ is a solution to (\ref{eq:intsys}) for some $T>0$. Then
\begin{align}
\sup_{(x,t)\in K_T}&|u(x,t)| \label{eq:basiclemma}\\ &
\leq \Big(\sup_{x\in K_0}|u_0(x)| + T\sup_{(x,t)\in K_T}|g(x,t)|\Big)\exp\big(T\sup_{(x,t)\in K_T}|f(x,t)|\big).\nonumber
\end{align}
\end{lemma}
\emph{Proof.} Apply Gronwall's inequality to the function
$s \to \sup_{(x,t)\in K_s}|u(x,t)|$. \hfill $\Box$

\begin{rem}\label{rem:Picard}
As preparation for the random case it is important to note that the solution to \eqref{eq:intsys} can be obtained by Picard iteration, which converges in $\cC(K_T)$ for arbitrarily chosen base interval $K_0$ and time $T$. To see this, we write equation \eqref{eq:intsys} in somewhat abusive, but suggestive notation as
\[
u(x,t) = u_0\big\vert_{\gamma(x,t,0)} + \int_0^t\big(fu+g\big)\big\vert_{(\gamma(x,t,\tau),\tau)}\dd\tau
\]
The Picard iteration scheme is
\begin{eqnarray} \label{eq:iterate}
u^{(0)}(x,t) & = & u_0\big\vert_{\gamma(x,t,0)}\nonumber\\
u^{(n+1)}(x,t) & = & u_0\big\vert_{\gamma(x,t,0)} + \int_0^t\big(fu^{(n)}+g\big)\big\vert_{(\gamma(x,t,\tau),\tau)}\dd\tau
\end{eqnarray}
Introducing the bounds
\[
   C = \sup_{(x,t)\in K_T}|f(x,t)|,\quad
   M = \sup_{(x,t)\in K_T}|g(x,t)| + C\sup_{x\in K_0}|u_0(x)|
\]
one gets the estimates
\[
   \sup_{(y,s)\in K_t}|u^{(1)}(y,s) - u^{(0)}(y,s)|\leq Mt\quad \mbox{for}\quad 0\leq t \leq T
\]
and, by induction,
\[
   \sup_{(y,s)\in K_t}|u^{(n+1)}(y,s) - u^{(n)}(y,s)|\leq M C^n t^{n+1}/(n+1)!
   \quad \mbox{for}\quad 0\leq t \leq T.
\]
Indeed, the induction step follows from the estimate
\[
\sup_{(y,s)\in K_t}|u^{(n+1)}(y,s) - u^{(n)}(y,s)|\leq
C \int_0^t \sup_{(y,s)\in K_\tau}|u^{(n)}(y,s) - u^{(n-1)}(y,s)|\dd\tau,
\]
which is based on the geometry of the characteristic curves. It follows that the sequence
\[
   u^{(n)}(x,t) = u^{(0)}(x,t) + \sum_{j=1}^n\big(u^{(j)}(x,t) - u^{(j-1)}(x,t)\big)
\]
converges absolutely on $K_T$. Passing to the limit in \eqref{eq:iterate} one sees that the limit coincides with the unique solution to \eqref{eq:intsys}. More details on these type of arguments can be found in \cite[Section 10]{Petrovsky}.
\end{rem}

%
%
\section{Random classical solutions}
\label{Sec:classrand}
%
%
Let $(\Omega,\Sigma,P)$ be a probability space. A \emph{random field} $A$ on $\R^d$ is a map
\begin{equation}\label{eq:rf}
   \R^d\times\Omega \to \R, (x,\omega) \to A(x,\omega)
\end{equation}
such that the map $\omega \to A(x,\omega)$ is measurable for every $x  \in \R^d$, more precisely, measurable with respect to $\Sigma$ and the Borel $\sigma$-algebra $\cB(\R)$.
The paths or trajectories of the random field are the maps $x \to A(x,\omega)$.

We may consider random fields with continuous or differentiable paths.
It is well known that a random field $A$ with almost surely continuous paths is jointly measurable, that is, the map in \eqref{eq:rf} is measurable with respect to $\cB(\R^d)\times\Sigma$.

\begin{rem}\label{rem:RF}
Suppose $A$ is a random field on $\R^d$ and $B$ a random field on $\R$, both with almost surely continuous paths. Then one may define the composed random field
\[
    (x,\omega) \to B\big(A(x,\omega),\omega\big).
\]
As the composition of the Borel measurable maps $(x,\omega) \to \big(A(x,\omega),\omega\big)$ and $(y,\eta) \to B(y,\eta)$ it is also jointly measurable (and it has almost surely continuous paths).
\end{rem}

We apply these observations to system \eqref{eq:hsys} with random coefficients. We make the assumption that the entries of $\lambda, f, g$ are random fields on $\R^2$ with almost surely continuous paths, and similarly for the initial data. In addition, the entries $\lambda_i(x,t,\omega)$ are assumed to satisfy condition (LB), with additional care given to the dependence on $\omega$. We require the following conditions:
\begin{enumerate}
\item[(AL)] (Almost sure local Lipschitz property)
For every compact subset $K\subset \R$, every compact time interval $I$, and almost all $\omega\in\Omega$, there is $\ell(\omega)>0$ such that, for all $i = 1,\ldots, n$,
\[
   |\lambda_i(x,t,\omega) - \lambda_i(y,t,\omega)| \leq \ell(\omega)|x-y|\quad\mbox{for\ all\ } x,y\in K, t\in I.
\]
\item[(UB)] (Uniform global bound)
For every compact time interval $I$ there is a constant $c>0$ such that, for all $i = 1,\ldots, n$ and for almost all $\omega\in\Omega$,
\[
   |\lambda_i(x,t,\omega)| \leq c \quad\mbox{for\ all\ } x\in \R, t\in I.
\]
\end{enumerate}
The integral operator arising in solving the $i$th characteristic equation now reads
\[
   (\cM h)(t,\omega) = x_0 -\int_t^{t_0} \lambda_i(h(s,\omega),s,\omega)\dd s.
\]
If $h(s,\omega)$ is a jointly measurable process with continuous paths, the same holds for $\lambda_i(h(s,\omega),s,\omega)$, using Remark\;\ref{rem:RF}. These properties are also preserved after integration, so that
$(\cM h)(t,\omega)$ will have the same property. Thus the integral operator preserves the measurability. Similarly to the procedure outlined in Remark\;\ref{rem:Picard}, the solution
$h(t,\omega) = \gamma_i(x_0,t_0,t,\omega)$ at fixed $(x_0,t_0)$ can be obtained by successive iteration
\[
   h^{(0)} = x_0,\quad h^{(n+1)} = x_0 + \int_t^{t_0} \lambda_i(h^{(n)}(s,\omega),s,\omega)\dd s
\]
which converges in $\cC(I)$ for any compact time interval around $t_0$, which can be chosen independently of $\omega$. This proves that $\omega\to\gamma_i(x_0,t_0,t,\omega)$ is measurable.
Continuity/differentiability with respect to $(x_0,t_0,t)$ at fixed $\omega$ follow as in Lemma\;\ref{lemma:char}.

We turn to the random version of system \eqref{eq:hsys} and \eqref{eq:intsys}. All participating functions are now random fields, that is, $u = u(x,t,\omega)$ etc. Thanks to condition (UB), the domains of determinacy \eqref{eq:KT} do not depend on $\omega$.
\begin{prop}
\label{prop:randsolu}
Let $\lambda(x,t,\omega)$, $u_0(x,\omega)$, $f(x,t,\omega)$ and $g(x,t,\omega)$ be random fields with almost surely continuous paths on $\R$ and on $\R^2$, respectively. Further, assume that $\lambda$ satisfies conditions (AL) and (UB).
Then there is an almost surely unique process $u(x,t,\omega)$ on $\R^2$ with continuous paths which solves system \eqref{eq:intsys} for almost all $\omega\in\Omega$.

In addition, if the paths of $\lambda$, $u_0$, $f$ and $g$ are $k$-times continuously differentiable ($k\geq 1$), so are the paths of $u$, and the process $u(x,t,\omega)$ is a pathwise solution of the system of differential equations \eqref{eq:hsys}.
\end{prop}
\emph{Proof.} The proof follows the lines of Proposition\;\ref{prop:classsolu}, making use of Remark\;\ref{rem:Picard}. One may fix a domain of determinacy $K_T$ and perform the iteration scheme
\eqref{eq:iterate}, where $\gamma$, $u_0$, $f$, $g$, and $u^{(n)}$ now depend on $\omega$. Using Remark\;\ref{rem:RF}, the integrands on the right-hand side of \eqref{eq:iterate} are continuous, measurable random fields, hence so is $u^{(n+1)}$. Uniform convergence of the iteration scheme implies that the limit $u(x,t,\omega)$ is a jointly measurable random field with continuous paths, which solves \eqref{eq:intsys}. Since the choice of $K_T$ is arbitrary, the solution is global. Uniqueness and differentiability properties are obtained exactly as in Proposition\;\ref{prop:classsolu}. \qed

Our next goal will be obtaining estimates on the moments of the solution and its derivatives, provided the data have finite moments. These estimates will be at the core of the proof of existence of a generalized solution in Section\;\ref{Sec:exis}, so we state them in some detail.

Whenever needed, we assume from now on that the probability space $(\Omega,\Sigma,P)$ is complete and that separable versions of the involved processes are chosen. If $u$ is a separable process with continuous paths and $K$ is the closure of a bounded open subset of $\R^2$, then $\sup_{(x,t)\in K}|u(x,t,\omega)|$ is a measurable function on $\Omega$.
\begin{lemma}\label{lemma:basicrandlemma}
Let $u(x,t,\omega)$ be the process solving \eqref{eq:intsys} constructed in Proposition\;\ref{prop:randsolu} under the hypotheses stated there. Let $T > 0$ and $K_T$ be the trapezoidal domain as in \eqref{eq:KT} bounded by the lines of slope $\pm c$ with $c$ from condition (UB). Then:

(a) The estimate \eqref{eq:basiclemma} of Lemma\;\ref{lemma:basiclemma} holds pathwise almost surely.

(b) For $1 \leq p < \infty$, the $p$-th moments satisfy
\begin{align}
\rE\Big(&\sup_{(x,t)\in K_T}|u(x,t)|^p\Big) \label{eq:pthmoment}\\ &
     \leq 2^p\, \rE\Big(\big(\sup_{x\in K_0}|u_0(x)|^p + T^p\sup_{(x,t)\in K_T}|g(x,t)|^p\big)\exp\big(pT\sup_{(x,t)\in K_T} |f(x,t)|\big)\Big). \nonumber
\end{align}
\end{lemma}
\emph{Proof.} Part (a) is clear. Part (b) is obtained by taking $p$th powers in inequality \eqref{eq:basiclemma} and integrating over $\Omega$. \qed
%
%
\section{Colombeau random fields}
\label{Sec:Colombeaurandfun}
%
%
In this section we recall basic facts from the framework of Colombeau algebras of generalized functions \cite{c1, c2}.
The basic reference is \cite{GKOS}; the spaces and algebras of Colombeau random fields (also termed Colombeau stochastic processes)
are modelled after \cite{GOPS2018}. We shall derive an important new characterization and also prove that the moments of certain Colombeau stochastic processes can be defined as Colombeau generalized functions.
While Colombeau algebras can be defined on any open subset of $\R^d$ -- even on any smooth manifold with or without boundary -- it will suffice for our purposes to consider Colombeau algebras defined on $\R^d$ only.

We begin with the deterministic case. A family of real or complex numbers $(r_\eps)_{\eps \in (0,1]}$ is called \emph{moderate}
if
\[
   \exists a\geq 0: |r_\eps| = \cO(\eps^{-a})\ \rm{as}\ \eps \to 0.
\]
It is called \emph{negligible} if
\[
   \forall b\geq 0: |r_\eps| = \cO(\eps^{b})\ \rm{as}\ \eps \to 0.
\]
The sets of moderate and negligible families of numbers are denoted by $\cE_M$ and $\cN$, respectively.

The basic objects of the theory as we use it are families
$(u_\eps)_{\eps \in (0,1]}$ of smooth functions $u_\eps \in \Cinf(\R^d)$ for
$0 < \eps \leq 1$.
We single out the following subalgebras:

{\em Moderate families of smooth functions}, denoted by $\EM(\R^d)$, are defined by the property:
\begin{equation}
  \forall K \Subset \R^d\,\forall \alpha \in \N_0^d:
    \;\big(\textstyle\sup_{x\in K} |\p^\alpha u_\eps(x)|\big)_{\eps \in (0,1]} \in \cE_M
           \label{eq:mofu}
\end{equation}
{\em Negligible families of smooth functions}, denoted by $\cN(\R^d)$, are defined by the property:
\begin{equation}
  \forall K \Subset \R^d\,\forall \alpha \in \N_0^d:
    \;\big(\textstyle\sup_{x\in K} |\p^\alpha u_\eps(x)|\big)_{\eps \in (0,1]} \in \cN
             \label{eq:nufu}
\end{equation}
In other words, moderate families satisfy a locally uniform polynomial estimate in $(\eps^{-1})$ as $\eps \to 0$,
together with all derivatives, while negligible families
vanish faster than any power of $\eps$ in the same situation. The negligible families
form a differential ideal in the collection of moderate families. The {\em Colombeau
algebra of generalized functions} is the factor algebra
\[
   \cG(\R^d) = \EM(\R^d)/\cN(\R^d).
\]
The algebra $\cG(\R^d)$ just defined coincides with the {\em special Colombeau
algebra} in \cite[Def. 1.2.2]{GKOS}, where the notation $\cG^s(\R^d)$
has been employed. It was called the {\em simplified Colombeau algebra} in \cite{Biag}.

In the zero-dimensional case the factor algebra $\cE_M/\cN$ is called the \emph{ring of Colombeau generalized numbers} and usually denoted by $\widetilde{\R}$
or $\widetilde{\C}$. It coincides with the ring of constants in the algebras $\cG(\R^d)$.

\emph{Restriction to hyperplanes.} The restriction of an element $u \in \cG(\R^d)$ to a hyperplane $H$, e.g. given by $x_d = 0$
is defined on representatives by
\[
   u\vert_H = \ \mbox{class of}\ \big(u_\eps(\cdot,0)\big)_{\eps \in (0,1]}.
\]
\emph{Embedding of distributions.} The space of Schwartz distributions is embedded in $\cG(\R^d)$ by convolution:
\begin{equation} \label{eq:imbedding}
   \iota:\cD'(\R^d) \to \cG(\R^d),\;
     \iota(v) = \ \mbox{class of}\ \big(v \ast (\chi\rho_\eps)\big)_{\eps \in (0,1]},
\end{equation}
where $\rho \in \cS(\R^d)$ is a fixed test function with $\int\rho(x)\dd x = 1$, $\int x^\alpha\rho(x)\dd x = 0$ for all $\alpha\in\N_0^d$, $|\alpha|\geq 1$, and
\begin{equation} \label{eq:molli}
   \rho_\eps(x) = \eps^{-d}\rho\left(x/\eps\right).
\end{equation}
Further, $\chi \in \cD(\R^d)$ is a cut-off function which is identically equal to one in a neighborhood of zero.
The embedding respects derivatives and renders $\C^\infty(\R^d)$ a faithful subalgebra. More precisely, let
\begin{equation} \label{eq:standardimbedding}
   \sigma:\cC^\infty(\R^d) \to \cG(\R^d),\;
     \sigma(f) = \ \mbox{class of}\ [(\eps,x) \to f(x)]
\end{equation}
be the constant embedding. Then, for all $\alpha \in \N_0^d$,
\[
   \iota \circ \p^\alpha = \p^\alpha \circ \iota\qquad \mbox{and}\qquad \iota\vert_{\cC^\infty(\R^d)} = \sigma.
\]
Further, the embedding $\iota$ does not depend on the choice of the cut-off function $\chi$. If $v\in \cS'(\R^d)$, the cut-off $\chi$ can be omitted in formula \eqref{eq:imbedding}.
More details on the embedding can be found in the Appendix.

\emph{The association relation.} We end this excursion into the deterministic case by recalling the {\em association relation} on the Colombeau
algebra $\cG(\R^d)$. It identifies elements of $\cG(\R^d)$ if they
coincide in the weak limit. That is, $u_1, u_2 \in \cG(\R^d)$ are called \emph{associated},
$u_1 \approx u_2$, if
$ \lim_{\eps \to 0} \int\big(u_{1\eps}(x) - u_{2\eps}(x)\big) \psi(x)\,dx = 0 $
for all test functions $\psi \in \cD(\R^d)$. We shall also say that
$u$ is \emph{associated with a distribution} $v$ if $u_\eps \to v$
in the sense of distributions as $\eps\to 0$, that is, $u \approx \iota(v)$.

Turning to the stochastic case, let again $(\Omega,\Sigma,P)$ be a probability space. The vector and matrix norms are defined as in \eqref{eq:norms}.
For $1\leq p < \infty$, the $L^p$-norm of a random variable, vector or matrix is
\[
    \|u\|_{L^p(\Omega)} = \Big(\rE\big(|u|^p\big)\Big)^{1/p} = \Big(\int_\Omega |u(\omega)|^p\dd P(\omega)\Big)^{1/p}
\]
and for $p = \infty$
\[
    \|u\|_{L^\infty(\Omega)} = \esssup_{\omega\in\Omega}|u(\omega)|.
\]
The space $L^p(\Omega)$ comprises the measurable functions on $\Omega$ with finite $L^p$-norm. To unify notation, we will also write $L^0(\Omega)$ for the space of measurable functions.

The basic objects in the stochastic theory of Colombeau algebras are again families of random fields $(u_\eps)_{\eps \in (0,1]}$ on $\R^d$ with smooth paths.
Consider the following properties:
\begin{itemize}
\item[(S)] For almost all $\omega\in\Omega$ and all $0<\eps\leq 1$, the map $x\to u_\eps(x,\omega)$ belongs $\Cinf(\R^d)$.

\item[(L$_0$)] For all $x\in \R^d$ and $0<\eps\leq 1$, the map $\omega\to u_\eps(x,\omega)$ is measurable.
\end{itemize}
These two properties will be shared by all spaces to follow.
The spaces $\cE_{M,L^0}(\Omega,\R^d)$ and $\cN_{L^0}(\Omega,\R^d)$ are obtained by requiring the additional property
\begin{itemize}
\item[(M$_0$)] For almost all $\omega\in\Omega$, $\big(u_\eps(\cdot,\omega)\big)_{\eps\in (0,1]} \in \cE_M(\R^d)$
\end{itemize}
and, respectively,
\begin{itemize}
\item[(N$_0$)] For almost all $\omega\in\Omega$, $\big(u_\eps(\cdot,\omega)\big)_{\eps\in (0,1]} \in \cN(\R^d)$.
\end{itemize}
For $1\leq p \leq\infty$, properties (S) and (L$_0$) together with the further properties
\begin{itemize}
\item[(L$_{\rm p}$)] For all $x\in \R^d$, $\alpha \in \N_0^d$ and $0<\eps\leq 1$, the map $\omega\to \p^\alpha u_\eps(x,\omega)$ belongs to $L^p(\Omega)$ and, for all $K \Subset \R^d$,
$\sup_{x\in K} \|\p^\alpha u_\eps(x)\|_{L^p(\Omega)}$ is finite
\item[(M$_{\rm p}$)] For all $K \Subset \R^d$ and all $\alpha \in \N_0^d$, $\big(\sup_{x\in K} \|\p^\alpha u_\eps(x)\|_{L^p(\Omega)}\big)_{\eps \in (0,1]} \in \cE_M$
\item[(N$_{\rm p}$)] For all $K \Subset \R^d$ and all $\alpha \in \N_0^d$, $\big(\sup_{x\in K} \|\p^\alpha u_\eps(x)\|_{L^p(\Omega)}\big)_{\eps \in (0,1]} \in \cN$
\end{itemize}
define the spaces $\cE_{M,L^p}(\Omega,\R^d)$ and $\cN_{L^p}(\Omega,\R^d)$, respectively.
The spaces are referred to as \emph{moderate} and \emph{negligible} families of functions with values in $L^p(\Omega)$. The factor spaces
\[
    \cG_{L^p}(\Omega,\R^d) = \cE_{M,L^p}(\Omega,\R^d)/\cN_{L^p}(\Omega,\R^d)
\]
are the spaces of \emph{Colombeau stochastic processes} or \emph{Colombeau random fields with values in} $L^p(\Omega)$ ($p = 0$ or $1\leq p \leq\infty$).
The spaces are closed under differentiation; $\cG_{L^0}(\Omega,\R^d)$ and $\cG_{L^\infty}(\Omega,\R^d)$ are algebras, in addition.

Again, in the zero-dimensional case one obtains the spaces of \emph{Colombeau random variables} $\cG_{L^p}(\Omega) = \cE_{M,L^p}(\Omega)/\cN_{L^p}(\Omega)$.

\begin{rem}\label{rem:inclusions}
For $1 \leq p < p' \leq\infty$, the inclusions
\[
   \cE_{M,L^{p'}}(\Omega,\R^d) \subset \cE_{M,L^p}(\Omega,\R^d),\qquad \cN_{L^{p'}}(\Omega,\R^d)\subset \cN_{L^{p}}(\Omega,\R^d)
\]
are valid, but the corresponding maps on the factor spaces $\cG_{L^{p'}}(\Omega,\R^d) \to \cG_{L^{p}}(\Omega,\R^d)$ are not injective.

Further,
$\cE_{M,L^1}(\Omega,\R^d)$ is not contained in $\cE_{M,L^0}(\Omega,\R^d)$ so that there is no canonical map from $\cG_{L^1}(\Omega,\R^d)$ to $\cG_{L^0}(\Omega,\R^d)$.
This also implies that Colom\-beau random fields with values in $L^p(\Omega)$ ($p\geq 1)$ need not have paths which are Colombeau generalized functions; see
Example\;\ref{ex:counterexamples} in the Appendix for details.
\end{rem}

\emph{Embedding of generalized stochastic processes.} A \emph{random Schwartz distribution} is a weakly measurable map $v: \Omega\to \cD'(\R^d)$, i.e., for all $\psi\in\cD(\R^d)$, the maps $\omega\to \langle v(\omega),\psi\rangle$ are measurable. Performing the convolution in \eqref{eq:imbedding} pathwise, i.e., with $v(\omega)$ in place of $v$, produces an embedding of these processes into $\cG_{L^0}(\Omega,\R^d)$.
If in addition, $\langle v,\varphi\rangle \in L^p(\Omega)$ for all $\varphi\in \cD(\R^d)$ and the map $\varphi \to \langle v,\varphi\rangle$ is linear and continuous from $\cD(\R^d)$ to $L^p(\Omega)$, then $\iota(v)$
defines an element of $\cG_{L^p}(\Omega,\R^d)$. Again, details on these embeddings are worked out in the Appendix.

We proceed to an important observation which shows that in the defining properties (M$_{\rm p}$) and (N$_{\rm p}$), the supremum and the $L^p$-norm can be interchanged.
\begin{lemma}\label{lemma:supequiv}
Let $1\leq p \leq\infty$ and suppose a family $(u_\eps)_{\eps\in (0,1]}$ satisfies conditions (S) and (L$_{\rm p}$). Then condition (M$_{\rm p}$) is equivalent with
\begin{itemize}
\item[($\widetilde{M}_{p}$)] For all $K \Subset \R^d$ and all $\alpha \in \N_0^d$, $\big(\big\|\sup_{x\in K} |\p^\alpha u_\eps(x)|\big\|_{L^p(\Omega)}\big)_{\eps \in (0,1]} \in \cE_M$.
\end{itemize}
The same statement holds for (N$_{\rm p}$) and the corresponding condition ($\widetilde{N}_{\rm p}$).
\end{lemma}
\emph{Proof.} We show the argument only for $|\alpha| = 0$. Consider first the case $1\leq p < \infty$ and let $u$ be a random field on $\R^d$ with smooth paths
such that, for all $K \Subset \R^d$, $\sup_{x\in K} \|\p^\alpha u(x)\|_{L^p(\Omega)}$ is finite.

If $K$ is a cuboid of the form $K = [a_1,b_1] \times \dots \times [a_d,b_d]$ then
\[
   u(x_1, \ldots, x_d,\omega) = u(a_1, \ldots, a_d,\omega) + \int_{a_1}^{x_1}\ldots \int_{a_d}^{x_d} \p^\beta u(y_1, \ldots, y_d,\omega)\dd y_1\ldots \dd y_d
\]
with $\beta = (1,\ldots,1)$. This implies (writing $a = (a_1,\ldots, a_d)$, $x = (x_1,\ldots, x_d)$, $y = (y_1,\ldots, y_d)$) that
\begin{equation}\label{eq:supubound}
   \sup_{x\in K}|u(x,\omega)| \leq |u(a,\omega)| + \int_K |\p^\beta u(y,\omega)|\dd y
\end{equation}
and
\begin{align*}
 \big\|\sup_{x\in K}& |u(x)|\big\|_{L^p(\Omega)} = \left(\int_\Omega \big(\sup_{x\in K}|u(x,\omega)|\big)^p\dd P(\omega)\right)^{1/p} \\
   &\leq \|u(a,\cdot)\|_{L^p(\Omega)} + \left(\int_\Omega \left(\int_K |\p^\beta u(y,\omega)|\dd y\right)^p \dd P(\omega)\right)^{1/p}\\
   &\leq \|u(a,\cdot)\|_{L^p(\Omega)} + \int_K\left(\int_\Omega |\p^\beta u(y,\omega)|^p \dd P(\omega)\right)^{1/p}\dd y\\
   &\leq \|u(a,\cdot)\|_{L^p(\Omega)} + \int_K\sup_{x\in K}\left(\int_\Omega |\p^\beta u(x,\omega)|^p \dd P(\omega)\right)^{1/p}\dd y\\
   &\leq \|u(a,\cdot)\|_{L^p(\Omega)} + \nu(K) \sup_{x\in K}\|\p^\beta u(x,\cdot)\|_{L^p(\Omega)}
\end{align*}
where $\nu$ denotes Lebesgue measure; to obtain the third line, we used Minkow\-ski's inequality for integrals.
An arbitrary compact subset $K\subset \R^d$ can be covered by finitely many cuboids as above, so we get the same inequality with a possibly larger compact set on the right-hand side.
The reverse inequality
\[
   \sup_{x\in K}\|u(x,\cdot)\|_{L^p(\Omega)} \leq \big\|\sup_{x\in K} |u(x)|\big\|_{L^p(\Omega)}
\]
is always true. Consequently, if $u$ is replaced by a family $(u_\eps)_{\eps\in (0,1]}$, moderate bounds in (M$_{\rm p}$) imply those in ($\widetilde{M}_{\rm p}$), and vice versa, and the same holds for (N$_{\rm p}$) and ($\widetilde{N}_{\rm p}$).
The case $p = \infty$ follows from an obvious modification of the norm estimates.
\qed

We turn to proving that the $p$th moments are well-defined as Colombeau functions for elements of $\cG_{L^p}(\Omega,\R^d)$. Clearly, taking any representative $(u_\eps)_{\eps\in (0,1]}$ and $x\in \R^d$,
$\|u_\eps(x,\cdot)\|_{L^p(\Omega)}$ is finite. The question is whether $(\eps,x)\to \rE(u_\eps(x,\cdot)^p)$ defines an element of $\cG(\R^d)$.
\begin{prop}\label{prop:moments}
Let $1 \leq p'\leq \infty$ and $u \in \cG_{L^{p'}}(\Omega,\R^d)$. Then $\rE(u^p)$ is a well-defined element of $\cG(\R^d)$ for every integer $p \leq p'$.
\end{prop}
\emph{Proof.} Pick a representative $(u_\eps)_{\eps\in(0,1]} \in \cE_{M,L^p}(\Omega,\R^d)$ of $u$. In a first step, we have to show that the map $x \to \rE(u_\eps^p(x))$ is smooth for fixed $\eps$.
For whatever multi-index $\alpha$, $\p^\alpha u_\eps^p$ is a sum of terms each of which contains exactly $p$ factors $v_{1,\eps}, \ldots, v_{p,\eps}$ (counted with multiplicity), where each $v_{j,\eps}$ is some derivative of $u_\eps$ of order less or equal to $|\alpha|$. Each of those factors belongs to $L^{p'}(\Omega)$, hence to $L^{p}(\Omega)$.
Thus the generalized H\"{o}lder inequality with $p$ factors\footnote{If $f_1,\ldots, f_n$ are in $L_{p_j}$ with $\sum_{j=1}^n\frac{1}{p_j} = 1$, then
$\|\prod_{j=1}^n f_j\|_{L_1} \leq \prod_{j=1}^n\|f_j\|_{L_{p_j}}$.
Apply this with $n = p$ and $p_j = p$ for all $j$.}
implies that $\p^\alpha u_\eps^p$ belongs to $L^1(\Omega)$. Let $K$ be a compact subset of $\R^d$. Then
\[
  \sup_{x\in K}\|v_{1,\eps}\cdots v_{p,\eps}\|_{L^1(\Omega)} \leq \sup_{x\in K}\|v_{1,\eps}\|_{L^p(\Omega)}\cdots \sup_{x\in K}\|v_{p,\eps}\|_{L^p(\Omega)},
\]
which in turn shows that $(\p^\alpha u_\eps^p)_{\eps\in(0,1]}$ is an element of $\cE_{M,L^1}(\Omega, \R^d)$, i.e.,
\begin{equation}\label{eq:L1moderate}
\big(\sup_{x\in K}\|\p^\alpha u_\eps^p(x)\|_{L^1(\Omega)}\big)_{\eps\in(0,1]} \in \cE_M.
\end{equation}
Invoking Lemma\;\ref{lemma:supequiv} we have that
\begin{equation*}
\big(\big\|\sup_{x\in K} |\p^\alpha u_\eps^p(x)|\big\|_{L^1(\Omega)}\big)_{\eps \in (0,1]} \in \cE_M.
\end{equation*}
In particular, $\sup_{x\in K} |\p^\alpha u_\eps^p(x)|$ is integrable so that interchanging integration and differentiation in the following computation is justified:
\[
  \p^\alpha\rE(u_\eps^p(x)) = \p^\alpha\int_\Omega u_\eps^p(x,\omega)\dd P(\omega) = \int_\Omega\p^\alpha u_\eps^p(x,\omega)\dd P(\omega).
\]
Thus the function $x\to \rE(u_\eps^p(x))$ is smooth. At the same time, equation \eqref{eq:L1moderate} shows that the family $\big(\rE(u_\eps^p)\big)_{\eps\in(0,1]}$
is moderate and hence defines an element $\rE(u^p)$ in $\cG(\R^d)$.

Finally, it remains to show that this element does not depend on the choice of representative. Thus let $(v_\eps)_{\eps\in(0,1]} \in \cN_{L^p}(\Omega,\R^d)$ be a negligible family. We have to show that the family $\rE\big((u_\eps + v_\eps)^p\big) - \rE\big(u_\eps^p\big)$ defines an element of $\cN(\R^d)$. Let $K$ be a compact subset of $\R^d$. Using the inequality $|(a+b)^p - a^p| \leq p(|a| + |b|)^{p-1}|b|$ leads to the estimate
\begin{align*}
\sup_{x\in K}\big|\rE\big((u_\eps + v_\eps)^p\big) & - \rE\big(u_\eps^p\big)\big| \\
         &\leq \sup_{x\in K} \int_\Omega p\big(|u_\eps(x,\omega)| + |v_\eps(x,\omega)|\big)^{p-1}|v_\eps(x,\omega)|\dd P(\omega)\nonumber\\
         &\leq \sup_{x\in K} p\big\||u_\eps(x)| + |v_\eps(x)|\big\|_{L^p(\Omega)}^{1/q} \|v_\eps(x)\|_{L^p(\Omega)} \nonumber
\end{align*}
where $\frac1{p} + \frac1{q} = 1$; to obtain the last line we used H\"{o}lder's inequality, observing that $q = \frac{p}{p-1}$. The right-hand side is the product of a moderate and a negligible term, hence negligible. By Theorem 1.2.3 of \cite{GKOS}, the zero order estimate suffices to infer the negligibility also for the derivatives. \qed
\begin{lemma}\label{lemma:tensor}
The map $\pi:u\to u\otimes u$ from $\cG_{L^2}(\Omega,\R^d)$ to $\cG_{L^1}(\Omega,\R^d\times \R^d)$ is injective.
\end{lemma}
\emph{Proof.} Let $K$ be a compact subset of $\R^d$ and $\alpha,\beta$ multi-indices. The inequality
\[
    \sup_{(x,y)\in K\times K}\|\p^\alpha_x\p^\beta_y u_\eps(x)u_\eps(y)\|_{L^1(\Omega)}
       \leq \sup_{x\in K}\|\p^\alpha_x u_\eps(x)\|_{L^2(\Omega)}\sup_{x\in K}\|\p^\beta_y u_\eps(y)\|_{L^2(\Omega)}
\]
shows that, on the level of representatives, $\pi$ maps $\cE_{M,L^2}(\Omega,\R^d)$ and $\cN_{L^2}(\Omega,\R^d)$ into $\cE_{M,L^1}(\Omega,\R^d\times \R^d)$ and
$\cN_{L^1}(\Omega,\R^d\times \R^d)$, respectively. It remains to show that if $(u_\eps)_{\eps\in (0,1]} \in \cE_{M,L^2}(\Omega,\R^d)$ and
$(u_\eps\otimes u_\eps)_{\eps\in (0,1]} \in \cN_{L^1}(\Omega,\R^d\times \R^d)$ then
$(u_\eps)_{\eps\in (0,1]} \in \cN_{L^2}(\Omega,\R^d)$. But this follows from the equality
\[
   \big\| \sup_{x\in K} |\p^\alpha_x u_\eps(x)| \big\|^2_{L^2(\Omega)} = \big\| \sup_{x\in K} \sup_{y\in K}|\p^\alpha_x u_\eps(x)||\p^\alpha_y u_\eps(y)| \big\|_{L^1(\Omega)}
\]
and Lemma\;\ref{lemma:supequiv}. \qed
\begin{cor}\label{cor:ACOV}
Let $u \in \cG_{L^2}(\Omega,\R^d)$. Then the autocovariance function
\begin{equation}\label{eq:ACOV}
   (x,y)\to C(x,y) = \rE\big(u(x)u(y)\big) - \rE(u(x))\rE(u(y))
\end{equation}
is a well-defined element of $\cG(\R^d\times \R^d)$.
\end{cor}
\emph{Proof.} By Lemma\;\ref{lemma:tensor}, $u\otimes u$ belongs to $\cG_{L^1}(\Omega,\R^d\times \R^d)$. The assertion follows from Proposition\;\ref{prop:moments}.
 \qed
%
%
\section{Existence/uniqueness of generalized solutions}
\label{Sec:exis}
%
%
In this section we will construct Colombeau stochastic processes in $\cG_{L^p}(\Omega,\R^2)$ which solve the hyperbolic system \eqref{eq:hsys}. The value of $p$ is
either $p=0$ or $1\leq p \leq\infty$ throughout. It will be useful to single out various types of processes that will enter as coefficients in the system.

\begin{defn}\label{def:bounded}
(a) A stochastic process $v\in \cG_{L^p}(\Omega,\R^2)$ is called \emph{uniformly bounded on strips} if it has a representative with the property:
\begin{itemize}
\item[{}] For every compact time interval $I$ there is a constant $c>0$ such that, for almost all $\omega\in\Omega$ and all $\eps\in (0,1]$,
\begin{equation}\label{eq:globalbound}
\sup_{x\in\R, t\in I}|v_\eps(x,t,\omega)| \leq c.
\end{equation}
\end{itemize}
(b) A stochastic process $v\in \cG_{L^p}(\Omega,\R^2)$ is \emph{of uniform local logarithmic type} if it has a representative with the property:
\begin{itemize}
\item[{}] For every compact subset $K\subset \R^2$ there is a constant $c>0$ such that, for almost all $\omega\in\Omega$ and all $\eps\in (0,1]$,
\begin{equation}\label{eq:uniflogtype}
\sup_{(x,t)\in K}|v_\eps(x,t,\omega)| \leq c|\log\eps|.
\end{equation}
\end{itemize}
(c) A stochastic process $v\in \cG_{L^p}(\Omega,\R^2)$ is \emph{of almost certain local logarithmic type} if it has a representative with the property:
\begin{itemize}
\item[{}] For every compact subset $K\subset \R^2$ and almost all $\omega\in\Omega$ there is $c(\omega)>0$ such that, for all $\eps\in (0,1]$,
\begin{equation}\label{eq:loclogtype}
\sup_{(x,t)\in K}|v_\eps(x,t,\omega)| \leq c(\omega)|\log\eps|.
\end{equation}
\end{itemize}
(d) A stochastic process $v\in \cG_{L^p}(\Omega,\R^2)$ is \emph{of almost certain local $L^1$-type} if it has a representative with the property:
\begin{itemize}
\item[{}] For every compact subset $K\subset \R^2$ and almost all $\omega\in\Omega$ there is $c(\omega)>0$ such that, for all $\eps\in (0,1]$,
\begin{equation}\label{eq:locL1type}
\iint_{K}|v_\eps(x,t,\omega)|\dd x \dd t \leq c(\omega).
\end{equation}
\end{itemize}
The definitions apply with obvious modifications to one-dimensional processes in $\cG_{L^p}(\Omega,\R)$.
\end{defn}
\begin{thm}\label{thm:gensolu0}
Assume that the entries of $\lambda$, $f$ and $g$ belong to $\cG_{L^0}(\Omega,\R^2)$, that $\lambda$ is uniformly bounded on strips and that $\p_x\lambda$ and $f$ are of almost certain local logarithmic type. Let $u_0 \in \cG_{L^0}(\Omega,\R)$. Then the Cauchy problem \eqref{eq:hsys} has a unique global solution $u\in \cG_{L^0}(\Omega,\R^2)$.
\end{thm}
\emph{Proof.} We choose representatives $u_{0\eps}, \lambda_\eps, f_\eps, g_\eps$ and write system \eqref{eq:hsys} in the form
	\begin{equation}
	\begin{array}{rcl}
	\left(\partial_t+\lambda_\eps(x,t,\omega)\partial_x\right)u_\eps(x,t,\omega) &=& f_\eps(x,t,\omega)u_\eps(x,t,\omega)+g_\eps(x,t,\omega)\\
	u_\eps(x,0,\omega)&=&u_{0\eps}(x,\omega)
    \end{array}
	\label{eq:hsysrep}
	\end{equation}
At fixed $\eps$, the coefficients are random fields with almost surely smooth paths. Because of its boundedness on strips and its local logarithmic type, $\lambda_\eps$ satisfies conditions (AL), (UB) of Section\;\ref{Sec:classrand} (with $c$ fixed and $\ell(\omega)$ possibly depending on $\eps$). By Proposition\;\ref{prop:randsolu} there is an almost surely unique process $u_\eps(x,t,\omega)$ with smooth paths which solves system \eqref{eq:hsysrep}.
This process will be the candidate as a representative of the generalized solution.

We begin by showing its moderateness. The assumption of global boundedness on strips guarantees that the slopes $\pm c$ of the bounding lines of the domains $K_T$, see \eqref{eq:KT}, can be chosen independently of $\omega$ and of $\eps$. It suffices to invoke Lemma\;\ref{lemma:basiclemma} to see that $\sup_{(x,t)\in K_T}|u_\eps(x,t,\omega)|$ is moderate for almost all $\omega$, using the moderate bounds on $u_{0\eps}$ and $g_\eps$ and the logarithmic bound on $f_\eps$.

We proceed by estimating the derivatives. For this purpose, it will be useful to collect some differential-algebraic relations. For clarity of exposition, we drop $x,t,\omega$ and $\eps$ in the formulas. Successively differentiating the system
\[
	\big(\partial_t+\lambda\partial_x\big)u = fu+g
\]
with respect to $x$ yields
\[
	\big(\partial_t+\lambda\partial_x\big)\p_xu = \big(f - (\p_x\lambda)\big)\p_x u + (\p_x f)u + \p_x g
\]
\[
	\big(\partial_t+\lambda\partial_x\big)\p_x^2u = \big(f - 2(\p_x\lambda)\big)\p_x^2 u + \big(2(\p_x f)- (\p_x^2\lambda)\big)\p_xu + (\p_x^2 f) u + \p_x^2 g
\]
We see by induction that the $k$th order $x$-derivative of $u_\eps(x,t,\omega)$ satisfies a linear hyperbolic system of the same form, with the highest order coefficient of the form
$(f_\eps - k\p_x\lambda_\eps)$ and the inhomogeneity being a linear function in $u_\eps$ and its $x$-derivatives up to order $k-1$. Hence repeated application of Lemma\;\ref{lemma:basiclemma} yields that $\sup_{(x,t)\in K_T}|\p_x^k u_\eps(x,t,\omega)|$ is moderate for almost all $\omega$.
Next, the $t$- and mixed derivatives can be simply estimated by differentiating the system of equations and using the already obtained estimates on the $x$-derivatives.
For example, up to order two the following differential-algebraic relations obtain:
\[
   \p_t u = -\lambda \p_x u + fu + g
\]
\[
  \p_t\p_x u =  -\lambda \p_x^2 + \big(f - (\p_x\lambda)\big)\p_x u + (\p_x f)u + \p_x g
\]
\[
  \p_t^2u = -\lambda \p_t\p_x u + (\p_t\lambda)\p_x u + f\p_t u + (\p_t f)u +\p_t g
\]
It is clear that one can go on step by step in a suitable order (for example, for $k = 3$ in the order $\p_t\p_x^2$, $\p_t^2\p_x$, $\p_t^3$) to obtain the moderate estimates on all derivatives of $u_\eps(x,t,\omega)$.

It remains to prove uniqueness. If $v = u_1 - u_2$ is the difference of two solutions $u_1, u_2$ in $\cG_{L^0}(\Omega,\R^2)$, then any representative of $v$ solves a system of the form
	\begin{equation}\label{eq:negsys}
	\begin{array}{rcl}
	\left(\partial_t+\lambda_\eps(x,t,\omega)\partial_x\right)v_\eps(x,t,\omega) &=& f_\eps(x,t,\omega)v_\eps(x,t,\omega)+n_\eps(x,t,\omega)\\
	v_\eps(x,0,\omega)&=&m_\eps(x,\omega)
    \end{array}
	\end{equation}
where the families $(n_\eps)_{\eps\in(0,1]}$ and $(m_\eps)_{\eps\in(0,1]}$ are negligible, i.e., they belong to $\cN_{L^0}(\Omega,\R^2)$ and
$\cN_{L^0}(\Omega,\R)$, respectively. Exactly the same strategy as in the proof of existence shows that $(v_\eps)_{\eps\in(0,1]}$ belongs to $\cN_{L^0}(\Omega,\R^2)$, that is, $u_1 = u_2$ in $\cG_{L^0}(\Omega,\R^2)$. \qed
\begin{thm}\label{thm:gensolup}
Let $1\leq p \leq \infty$. Assume that the entries of $\lambda$ and $f$ belong to $\cG_{L^\infty}(\Omega,\R^2)$, that $\lambda$ is uniformly bounded on strips and that $\p_x\lambda$ and $f$ are of uniform local logarithmic type. Assume also that $g\in \cG_{L^p}(\Omega,\R^2)$. Let $u_0 \in \cG_{L^p}(\Omega,\R)$. Then the Cauchy problem \eqref{eq:hsys} has a unique global solution $u\in \cG_{L^p}(\Omega,\R^2)$.
\end{thm}
\emph{Proof.}
We follow the same strategy as in the proof of Theorem\;\ref{thm:gensolu0}. The existence of a process $u_\eps(x,t,\omega)$ with smooth paths which solves \eqref{eq:hsysrep} was already  established there.

It remains to establish moderate bounds on
$\sup_{(x,t)\in K_T}\|\p^\alpha u_\eps(x,t)\big\|_{L^p(\Omega)}$ or, equivalently, on $\big\|\sup_{(x,t)\in K_T}|\p^\alpha u_\eps(x,t)|\big\|_{L^p(\Omega)}$ for all $\alpha \in \N_0^2$ (Lemma\;\ref{lemma:supequiv}).

For $\alpha = (0,0)$, these bounds are readily obtained from Lemma\;\ref{lemma:basicrandlemma}, using the  logarithmic estimate satisfied by
$\big\|\sup_{(x,t)\in K_T}|f_\eps(x,t,\omega)|\big\|_{L^\infty(\Omega)}$.

As in the proof of Theorem\;\ref{thm:gensolu0} we first estimate the $x$-derivatives inductively. The leading coefficient matrix in the equation for $\p_x^k u_\eps$ is $(f_\eps - k\p_x\lambda_\eps)$, which satisfies the required local logarithmic estimate. The lower order terms are multiplied by derivatives of $f_\eps$ and $\lambda_\eps$, which -- by assumption -- satisfy moderate estimates on $K_T$ in their $L^\infty(\Omega)$-norm. Hence the $L^p(\Omega)$-norms of the suprema over $K_T$ of the lower order terms are seen to be moderate from step number $(k-1)$. Finally, the estimates on the $t$- and mixed derivatives are successively obtained from the system of differential equations, using again the $L^\infty(\Omega)$-bounds on the derivatives of $f$ and $\lambda$.

To prove uniqueness, negligibility of $\big\|\sup_{(x,t)\in K_T}|\p^\alpha v_\eps(x,t)|\big\|_{L^p(\Omega)}$ is shown in the same way, based on system \eqref{eq:negsys}.
\qed
\begin{rem}\label{remark:only}
In specific cases the conditions on the coefficients can be relaxed; in particular, the logarithmic assumption can be avoided. We exemplify this in the simple case of the scalar transport equation
\begin{equation}\label{eq:trans}
(\p_t + \lambda(x,t)\p_x) u(x,t) = 0,\qquad u(x,0) = u_0(x)
\end{equation}
where $\lambda\in \cG_{L^0}(\Omega,\R^2)$ and $u_0\in \cG_{L^0}(\Omega,\R)$. Formally, the solution is given as the composition of $u_0$ with the base point of the characteristic curves, $u(x,t) = u_0(\gamma(x,t,0))$.
The question is whether this composition is a well-defined Colombeau function in $\cG_{L^0}(\Omega,\R^2)$. This is the case if $\gamma$ is \emph{c-bounded}, that is, for every compact subset $K\subset \R^2$ there is a compact subset $K'\subset\R$ such that $\gamma_\eps(x,t,0) \in K'$ for all $(x,t)\in K$ and $\eps\in (0,1]$, \cite[Proposition 1.2.8]{GKOS}. (In the stochastic case, $K'$ may depend on $\omega$). Alternatively, one can replace the Colombeau algebra $\cG$ by the
algebra $\cG^\tau$ of tempered Colombeau functions, in which setting compositions are always well-defined \cite[Proposition 1.2.29]{GKOS}.

Such cases arise when $\lambda$ depends on one of the variables $(x,t)$ only. For example, if $\lambda$ depends only on $x$ and there are constants $c_0, c_1>0$ such that for all $x\in \R$, almost all $\omega\in\Omega$ and all $\eps\in (0,1]$,
\begin{equation}\label{eq:upperlowerbound}
c_0 \leq |\lambda_{\eps}(x,\omega)| \leq c_1
\end{equation}
then $\gamma_\eps(x,t,0,\omega)$ is c-bounded independently of $\omega$ and \eqref{eq:trans} is uniquely solvable in $\cG_{L^0}(\Omega,\R^2)$. This follows by using the (deterministic) arguments given in the proof of \cite[Theorems 3.10]{DO2016}.
Switching to the $\cG^\tau$-setting, one can further relax the strict bounds in \eqref{eq:upperlowerbound} to $\eps^m \leq |\lambda_{\eps}(x,\omega)| \leq \eps^{-m}$ for some $m\in \N$ to construct a unique solution to \eqref{eq:trans} in $\cG_{L^0}^\tau(\Omega,\R^2)$, see the arguments in \cite[Section 2.2]{Schwarz:Thesis:2019}. Another specific instance where the special form of $\lambda$ can be exploited is in Example\;\ref{ex:wavegeometric}.

The case of a $t$-dependent coefficient is dealt with in the following proposition.
\end{rem}

\begin{prop}\label{prop:tonly}
Assume that in \eqref{eq:trans} the coefficient $\lambda = \lambda(t)$ depends on $t$ only, belongs to $\cG_{L^0}(\Omega,\R)$ and is of almost certain local $L^1$-type. Let $u_0\in \cG_{L^0}(\Omega,\R)$.
Then problem \eqref{eq:trans} has a unique solution $u\in \cG_{L^0}(\Omega,\R^2)$.
\end{prop}
\emph{Proof.}
The characteristic curves are given by
\[
   \gamma_{\eps}(x,t,\tau,\omega) = x - \int_\tau^{t}\lambda_{\eps}(s,\omega)\dd s.
\]
Since $\lambda$ is of almost certain local $L^1$-type, the random field $\gamma$ is almost surely c-bounded as a function of $x,t,\tau$. The discussion in Remark\;\ref{remark:only} shows that
\[
   u_\eps(x,t,\omega) = u_{0\eps}\big( \gamma_{\eps}(x,t,0,\omega)\big)
\]
is a well-defined representative of a solution in $u\in \cG_{L^0}(\Omega,\R^2)$. To prove uniqueness, let $v = u_1 - u_2$ be the difference of two solutions. Then
\[
   v_\eps(x,t,\omega) = m_\eps\big( \gamma_{\eps}(x,t,0,\omega)\big) + \int_0^t n_\eps\big(\gamma_{\eps}(x,t,\tau,\omega),\tau,\omega\big)\dd\tau
\]
where $m$ and $n$ are negligible. The desired negligibility of $v$ obviously follows from the c-boundedness of $\gamma_{\eps}(x,t,\tau,\omega)$.
\qed

We end this section by a statement of consistency of Colombeau solutions with classical solutions to \eqref{eq:hsys} when the latter exist. To fix notation, consider the classical version of \eqref{eq:hsys}
	\begin{equation}
	\begin{array}{rcl}
	\left(\partial_t+\lambda^\circ(x,t)\partial_x\right)u^\circ&=&f^\circ(x,t)u^\circ+g^\circ(x,t), \quad(x,t)\in\R^2,\\
	u^\circ(x,0)&=&u^\circ_0(x), \quad x\in\R,
	\end{array}
	\label{eq:hsysclass}
	\end{equation}
where $\lambda^\circ$, $f^\circ$, $g^\circ$, $u^\circ_0$ are random fields\footnote{From now on we use a circle as superscript to indicate classical processes; Colombeau processes will be denoted by small roman letters as before.} as compared with almost surely continuous paths. If $\lambda^\circ$ satisfies the conditions (AL) and (UB), then \eqref{eq:hsysclass} has an almost surely unique solution $u^\circ(x,t,\omega)$
with almost surely continuous paths (Proposition\;\ref{prop:randsolu}). Embedding $\lambda^\circ$, $f^\circ$, $g^\circ$ and $u^\circ_0$ into $\cG_{L^0}(\Omega,\R^2)$ and $\cG_{L^0}(\Omega,\R)$, respectively, by means of convolution with a truncated mollifier $\chi\rho_\eps$, see \eqref{eq:randomimbedding}, produces
elements $\lambda = \iota(\lambda^\circ)$, $\lambda_\eps = \lambda^\circ\ast(\chi\rho_\eps)$, and similarly for $f = \iota(f^\circ)$, $g = \iota(g^\circ)$, $u_0 = \iota(u^\circ_0)$. It follows from (UB) that $\lambda$ is uniformly bounded on strips. Further, the Lipschitz property (AL)
gives that $\p_x\lambda$ is locally bounded at fixed $\omega$ and thus is of almost certain local logarithmic type. The same applies to $f$ due the continuity of $f^\circ$. By Theorem\;\ref{thm:gensolu0} the Cauchy problem \eqref{eq:hsys} has a unique solution $u\in \cG_{L^0}(\Omega,\R^2)$. We are going to show that $u$ is associated (cf. Section\;\ref{Sec:Colombeaurandfun}) with $u^\circ$, at least when $u^\circ$ has $\cC^1$-paths. (By Proposition\;\ref{prop:randsolu} this is the case, if all coefficients have $\cC^1$-paths, in addition.)

\begin{prop}\label{prop:consistency}
Assume that the random fields  $\lambda^\circ$, $f^\circ$, $g^\circ$, $u^\circ_0$ have almost surely continuously differentiable paths and that $\lambda^\circ$ satisfies conditions (AL) and (UB). Then the generalized solution
$u\in \cG_{L^0}(\Omega,\R^2)$ is associated with the classical $\cC^1$-solution $u^\circ$. More precisely, the representatives $u_\eps$ of $u$ almost surely converge to $u^\circ$ in $\cC(\R^2)$.
\end{prop}
\emph{Proof.}
Subtracting the two systems
\[
   (\p_t + \lambda^\circ\p_x)u^\circ = f^\circ u^\circ + G,\qquad u^\circ(x,0) = u^\circ_0(x)
\]\[
   (\p_t + \lambda_\eps\p_x)u_\eps = f_\eps u_\eps + g_\eps,\qquad u_\eps(x,0) = u_{0\eps}(x)
\]
shows that the difference $u^\circ-u_\eps$ solves the system
\begin{eqnarray*}
  \p_t(u^\circ - u_\eps) &+& \lambda_\eps \p_x(u^\circ - u_\eps)\\[4pt]
     &=& f_\eps(u^\circ - u_\eps) + (\lambda_\eps - \lambda^\circ)\p_x u^\circ + (f^\circ - f_\eps) u^\circ + (g^\circ - g_\eps)
\end{eqnarray*}
with initial data $u^\circ_0 - u_{0\eps}$. The coefficient $\lambda_\eps$ has uniform bounds independently of $\eps$ and $\omega$. Thus Lemma\;\ref{lemma:basiclemma} shows that
\[
\begin{array}{l}
\|u^\circ-u_\eps \|_{L^\infty(K_T)}  \\[4pt]
 \quad\leq\Big(\|u^\circ - u_{0\eps}\|_{L^\infty(K_0)} + T\|(\lambda_\eps - \lambda^\circ)\p_x u^\circ + (f^\circ - f_\eps) u^\circ + (g^\circ - g_\eps)\|_{L^\infty(K_T)}\Big)\\[6pt]
  \qquad\qquad \cdot\exp\big(T\|f_\eps\|_{L^\infty(K_T)}\big).
\end{array}
\]
Lemma\;\ref{lem:uniformconvergence} below implies that the net $f_\eps$ is uniformly bounded on $K_T$ and all difference terms converge to zero, so the assertion follows.
\qed

One can relax the $\cC^1$-assumption on $u^\circ$ through a more careful analysis of the continuous dependence of the solution $u^\circ$ on $\lambda^\circ$ and the characteristic curves, as outlined, e.g., in 
\cite[Section 4.3]{Nedeljkovic:Thesis:2020}.

The following lemma turns out to be useful also at various places in Section\;\ref{Sec:apps}.

\begin{lemma}\label{lem:uniformconvergence}
Let $h\in\cC(\R^d)$ and $\chi$, $\rho_\eps$ as in \eqref{eq:imbedding}, \eqref{eq:molli}. Then $h\ast(\chi\rho_\eps) \to h$ locally uniformly as $\eps\to 0$.
\end{lemma}
\emph{Proof.} To be specific, assume that $\chi(x) \equiv 1$ for $|x|\leq a$ and $\chi(x) \equiv 0$ for $|x|\geq b$. We have
\begin{eqnarray*}
\big(h*(\chi\rho_\eps)\big)(x) - h(x) &=& \int h(x-y)\chi(y)\eps^{-d}\rho(y/\eps)\dd y - h(x)\\
&=& \int\big(h(x-\eps y)\chi(\eps y) - h(x)\big)\rho(y)\dd y. 
\end{eqnarray*}
We split the integral in an integration over $|y|\leq 1/\sqrt{\eps}$ and over $|y|\geq 1/\sqrt{\eps}$. 
In the first summand, $\chi(\eps y) \equiv 1$ when $\sqrt{\eps} \leq a$. If $x$ is taken from a compact set, then $x - \eps y$ remains in a slightly larger compact set for $|y|\leq 1/\sqrt{\eps}$, on which set $h$ is uniformly continuous.
Thus the first summand tends to zero. Recall that $\chi(\eps y) \equiv 0$ if $|y|\geq b/\eps$. Therefore, the second summand reads
\[
  \int_{1/\sqrt{\eps} \leq |y| \leq b/\eps}\big(h(x-\eps y)\chi(\eps y) - h(x)\big)\rho(y)\dd y.
\]
In the integration interval, $x-\eps y$ remains in a compact set as well. Thus $h(x-\eps y)\chi(\eps y) - h(x)$ is bounded, while $\int_{1/\sqrt{\eps}\leq |y|} |\rho(y)|\dd y$ tends to zero.
\qed

%
%
\section{Applications}
\label{Sec:apps}
%
%
\begin{ex}\label{ex:wavegeneral}
The one-dimensional wave equation with lower order terms, in its most general form, reads
\begin{equation}\label{eq:genwave}
\begin{array}{l}
   \p_t^2 u - \lambda^2(x,t)\p_x^2 u = f(x,t)u + h(x,t)\p_tu + k(x,t)\p_xu + g(x,t), \\[4pt] 
   u(x,0) = u_0(x), \quad \p_t u(x,0) = u_1(x),
\end{array}
\end{equation}
with $(x,t)\in\R^2$. We assume here that all coefficients as well as the initial data are Colombeau random fields. To place ourselves in the previous setting, we transform \eqref{eq:genwave} into a first order hyperbolic system, by introducing the new dependent variables
\[
   v = (\p_t - \lambda\p_x)u,\qquad w = (\p_t + \lambda\p_x)u.
\]
A simple calculation yields the equivalent system
\begin{equation}\label{eq:wavesyst}
\begin{array}{rcl}
  (\p_t + \lambda\p_x)v &=& \frac12\big(\frac{k}{\lambda} - \frac{\p_t \lambda}{\lambda} - \p_x\lambda\big)(w-v) + \frac{h}{2}(v+w) + fu + g,\\[5pt]
    (\p_t - \lambda\p_x)w &=& \frac12\big(\frac{k}{\lambda} + \frac{\p_t \lambda}{\lambda} - \p_x\lambda\big)(w-v) + \frac{h}{2}(v+w) + fu + g,\\[5pt]
    \p_t u &=& \frac12(v+w)
\end{array}
\end{equation}
with initial data
\[
\begin{array}{l}
v(x,0) = u_1(x) - \lambda(x,0)u_0'(x),\\[3pt]
w(x,0) = u_1(x) + \lambda(x,0)u_0'(x),\\[3pt]
u(x,0) = u_0(x).
\end{array}
\]
In order to apply Theorem\;\ref{thm:gensolu0} one has to assume that all coefficients belong to $\cG_{L^0}(\Omega, \R^2)$, $\lambda$ is uniformly bounded on strips and $\p_x\lambda$ as well as the combined coefficients on the right-hand side of
\eqref{eq:wavesyst} are of almost certain local logarithmic type. In case $k -\p_t\lambda\neq 0$ or $k +\p_t\lambda\neq 0$ one needs in addition that $\lambda$ is an invertible Colombeau function\footnote{A Colombeau function $\lambda$ is invertible, if it has a representative $(\lambda_\eps)_{\eps\in(0,1]}$ such that on every compact set $K$, it has a lower bound of the form $\inf_{y\in K}|\lambda_\eps(y)|\geq \eps^m$ for some $m\geq 0$ and sufficiently small $\eps$.} or process.
The conditions have to be modified in an obvious manner if one wants to apply Theorem\;\ref{thm:gensolup}.

We shall focus on a number of special, explicit examples.
\end{ex}
\begin{ex}\label{ex:wave}
The stochastic wave equation with random field coefficient depending on $x$ only:
\begin{equation}\label{stochasticWaveColombeau}
\begin{array}{l}
   \p_t^2 u - \lambda^2(x)\p_x^2 u = 0, \\[4pt]
   u(x,0) = u_0(x), \quad \p_t u(x,0) = u_1(x).
\end{array}
\end{equation}
In this case, the right-hand side of \eqref{eq:wavesyst} reduces to the sole term $\p_x\lambda(w-v)$ in the first two lines. Thus if $\lambda \in \cG_{L^0}(\Omega, \R)$, $\lambda$ is uniformly bounded and $\p_x\lambda$
of almost certain local logarithmic type, and the initial data belong to $\cG_{L^0}(\Omega, \R)$, the equation \eqref{stochasticWaveColombeau} has a unique solution $u \in \cG_{L^0}(\Omega, \R^2)$.

On the other hand, there also exist classical random field solutions, provided the coefficients and data are sufficiently regular; we exemplify this in the situation addressed in Proposition\;\ref{prop:randsolu}. To be specific, consider the classical equation
\begin{equation}\label{stochasticWaveClassical}
\begin{array}{l}
   \p_t^2 u^\circ - (\lambda^\circ)^2(x)\p_x^2 u^\circ = 0, \\[4pt]
   u^\circ(x,0) = u^\circ_0(x), \quad \p_t u^\circ(x,0) = u^\circ_1(x).
\end{array}
\end{equation}
Assume that $\lambda^\circ$, $u^\circ_0$ and $u^\circ_1$ are random fields on $\R$ with almost surely continuous paths such that condition (UB) is satisfied; assume further that the paths of $\lambda^\circ$ are continuously differentiable (then condition (AB) holds, in particular). Proposition\;\ref{prop:randsolu} guarantees the existence and almost sure uniqueness of a classical solution
to the integral equation corresponding to the first order system \eqref{eq:wavesyst}.

To relate the classical and the Colombeau solutions, consider $\lambda = \iota(\lambda^\circ)$, $u_0 = \iota(u^\circ_0)$, $u_1 = \iota(u^\circ_1)$, where $\iota$ is the embedding \eqref{eq:randomimbedding}.  It follows from (UB) that the Colombeau random process $\lambda$ is uniformly bounded (on strips). Because of the $\cC^1$-property of $\lambda^\circ$, the process $\p_x\lambda$ has an almost certain locally bounded representative (and in particular is of almost certain local logarithmic type). With the embedded coefficient and data, equation \eqref{stochasticWaveColombeau} has a unique solution $u \in \cG_{L^0}(\Omega, \R^2)$.

If in addition the paths of $\lambda^\circ$ and $u^\circ_0$ belong to $\cC^2(\R)$ and those of $u^\circ_1$ to $\cC^1(\R)$, then $u^\circ$ has paths in $\cC^1(\R^2)$, as can be read off from the corresponding system \eqref{eq:wavesyst}. Proposition\;\ref{prop:consistency} gives that $u$ is associated with $u^\circ$; actually $u_\eps(x,t,\omega) \to u^\circ(x,t,\omega)$ in
$\cC(\R^2)$ for almost all $\omega$ as $\eps\to 0$.
\end{ex}

\begin{ex}\label{ex:wavenoise}
The stochastic wave equation with constant coefficient and additive noise in the Colombeau sense:
\begin{equation}\label{eq:WaveAdditive}
\begin{array}{l}
   \p_t^2 u - \p_x^2 u = g(x,t), \\[4pt]
   u(x,0) = u_0(x), \quad \p_t u(x,0) = u_1(x).
\end{array}
\end{equation}
Assume that $u_0, u_1 \in \cG_{L^p}(\R)$ and $g\in \cG_{L^p}(\R^2)$ where $p=0$ or $1\leq p\leq\infty$. After reduction to system \eqref{eq:wavesyst}, Theorems\;\ref{thm:gensolu0} and \ref{thm:gensolup} immediately
show that \eqref{eq:WaveAdditive} has a unique solution u in $\cG_{L^p}(\R^2)$. Assume for simplicity that $u_0 = u_1 = 0$. Choosing representatives, it can be expressed, for $t\geq 0$, as
\[
  u_\eps(x,t) = \frac12\int_0^t\int_{x-t+s}^{x+t-s} g_\eps(y,s)\dd y\dd s = \frac12 \langle g_\eps, \eins_{\Gamma(x,t)}\rangle
\]
in the notation of distribution theory, where
\[
   \Gamma(x,t) = \{(y,s): 0 \leq s \leq t, x-t+s \leq y \leq x+t-s\}
\]
is the backward light cone. In particular, when $g$ is embedded space-time Gaussian white noise (Example\;\ref{ex:Brownian}(b)), $g = \iota(\dot{W})$, a representative is given by
\[
  u_\eps(x,t) =  \frac12 \langle \dot{W}\ast(\chi\rho_\eps), \eins_{\Gamma(x,t)}\rangle = \frac12 \langle \dot{W}, \eins_{\Gamma(x,t)}\ast(\check{\chi}\check{\rho}_\eps) \rangle
\]
where the check denotes inflection.
At fixed $(x,t)$, the function $(y,s)\to \eins_{\Gamma(x,t)}(y,s)$ is bounded and has compact support. A slight modification of Lemma\;\ref{lem:uniformconvergence} shows that
$\eins_{\Gamma(x,t)}\ast(\check{\chi}\check{\rho}_\eps)$ converges to $\eins_{\Gamma(x,t)}$ in $L^2(\R^2)$. The It\^o isometry \eqref{eq:ItoIsometry} gives that $u_\eps(x,t)$ is a Cauchy net in $L^2(\Omega)$
and thus it converges in mean square to a limit denoted by $\widehat{W}(x,t) = \langle \dot{W}, \eins_{\Gamma(x,t)}\rangle$. By a similar limiting argument, the covariance is given by
\[
   \rE(\widehat{W}(x_1,t_1)\widehat{W}(x_2,t_2)) = \nu(\Gamma(x_1,t_1) \cap \Gamma(x_2,t_2))
\]
where $\nu$ denotes Lebesgue measure.
This is well-known: $\widehat{W}$ is a modified Brownian sheet \cite[pp. 281 -- 283]{Walsh:1984}, and it is the unique random Schwartz distribution that is a pathwise weak solution of $\p_t^2 v - \p_x^2 v = \dot{W}$ and vanishes for $t<0$.
In summary, the Colombeau solution has the pointwise regularity property that, for every $(x,t)$ with $t> 0$, $u_\eps(x,t) \to \widehat{W}(x,t)$ in $L^2(\Omega)$.

\end{ex}

\begin{ex}\label{ex:wavegeometric}
The geometric wave equation: Let $C$ be the graph of a smooth function $y = c(x)$. Viewing $C$ as a one-dimensional Riemannian manifold, the wave equation on $C$ reads
\begin{equation}\label{eq:GeometricWave}
\begin{array}{l}
   \p_t^2 u - \lambda(x)\p_x (\lambda \p_x u) = 0, \\[4pt]
   u(x,0) = u_0(x), \quad \p_t u(x,0) = u_1(x)
\end{array}
\end{equation}
where
\begin{equation}\label{eq:transportcoeff}
\lambda(x) = (1 + c'(x)^2)^{-1/2}.
\end{equation}
Introducing $u^\pm = (\p_t \mp\lambda(x)\p_x) u$, the corresponding first order system decouples:
\begin{equation}\label{eq:geometricwavesyst}
\begin{array}{rcll}
  (\p_t + \lambda(x)\p_x)u^+ &=& 0,& u^+(x,0) = u_1(x) - \lambda(x)u_0'(x),\\[4pt]
    (\p_t - \lambda(x)\p_x)u^- &=& 0, & u^-(x,0) = u_1(x) + \lambda(x)u_0'(x),\\[4pt]
    \p_t u &=& \frac12(u^+ + u^-), & u(x,0) = u_0(x).
\end{array}
\end{equation}
One may add a stochastic perturbation to the curve $C$ and thus is led to considering $c\in\cG_{L^0}(\Omega,\R)$. Then $\lambda$ is automatically uniformly bounded.
It will follow from the analysis of the stochastic transport equation below that, due to the special form of $\lambda$, the logarithmic requirement of Theorem\;\ref{thm:gensolu0}
is not needed to infer existence and uniqueness of a solution:
Given initial data belonging to $\cG_{L^0}(\Omega,\R)$, system \eqref{eq:geometricwavesyst} has a unique solution in $\cG_{L^0}(\Omega,\R^2)$, and so has equation \eqref{eq:GeometricWave}.

{\bf Subcase, stochastic transport.} We start by analyzing the stochastic transport equation
\begin{equation}\label{eq:transportG}
 (\p_t + \lambda(x)\p_x)v = 0,\quad v(x,0) = v_{0}(x)
\end{equation}
with $\lambda$ of the form \eqref{eq:transportcoeff}. On the level of representatives, \eqref{eq:transportG} reads
\begin{equation}\label{eq:transport}
 (\p_t + \lambda_\eps(x,\omega)\p_x)v_\eps(x,t,\omega) = 0,\quad v_\eps(x,0,\omega) = v_{0\eps}(x,\omega)
\end{equation}
where $\lambda_\eps = (1 + c_\eps'(x,\omega)^2)^{-1/2}$ and $c_\eps$ is a representative of $c$. The case we have in mind is that $c = \iota(X)$ where $X$ is a random field belonging to $L^0(\Omega:\cD'(\R))$ so that 
$c_\eps = X \ast (\chi\rho_\eps)$. Note that classically much higher regularity of $c(x,\omega)$ is required to give a meaning to a transport equation with a coefficient of the form \eqref{eq:transportcoeff} \cite[Scction 4.2]{Flandoli:2011} or even to have a well-defined stochastic flow defining the characteristic curves \cite[Scction 4.4.1]{Flandoli:2011}.

The following calculations are all meant to be pathwise at fixed $\omega\in\Omega$; we drop $\omega$ in our notation for simplicity. The characteristic curves corresponding to \eqref{eq:transport} are obtained as solutions to
\[
  \frac{\dd}{\dd \tau} \gamma_\eps(x,t,\tau) = (1 + c_\eps'(\gamma_\eps(x,t,\tau))^2)^{-1/2},\quad \gamma_\eps(x,t,t) = x,
\]
and the unique solution to \eqref{eq:transport} is given by $v_{0\eps}(\gamma_\eps(x,t,0))$. Introducing $L_\eps(z) = \int_0^z \sqrt{1 + c_\eps'(y)^2}\dd y$, the relation
\begin{equation}\label{eq:arclength}
   L_\eps(x) - L_\eps(\gamma_\eps(x,t,0)) = t
\end{equation}
is readily obtained. At fixed $(x,t)$, the base point $\gamma_\eps(x,t,0)$ of the characteristic curve passing through $(x,t)$ is characterized by the property that the arclength of the curve $y= c_\eps(x)$ between $\gamma_\eps(x,t,0)$ and $x$ equals $t$.

At fixed $\eps$, the function $L_\eps$ is smooth and strictly monotonically increasing; we can compute
\[
   \gamma_\eps(x,t,0) = L_\eps^{-1}(L_\eps(x) - t).
\]
We are going to show that $\gamma_\eps(x,t,0)$ defines an element of $\cE_{M,L^0}(\Omega,\R^2)$ and is c-bounded at fixed $\omega$. Indeed,
\[
  \p_x\gamma_\eps(x,t,0) = \frac{L_\eps'(x)}{L_\eps'(L_\eps^{-1}(L_\eps(x) - t))} = \frac{\sqrt{1 + c_\eps'(x)^2}}{\sqrt{1 + c_\eps'(\gamma_\eps(x,t,0))^2}}
\]
and
\begin{equation}\label{eq:dtgamma}
  \p_t\gamma_\eps(x,t,0) =  \frac{(-1)}{\sqrt{1 + c_\eps'(\gamma_\eps(x,t,0))^2}}.
\end{equation}
The denominators are uniformly bounded, as are the denominators arising after higher order differentiation. Thus $\p_x\gamma_\eps(x,t,0)$ and $\p_t\gamma_\eps(x,t,0)$
are moderate. Further,
\[
   \gamma_\eps(x,t,0)  =  x + \int_0^t \p_t\gamma_\eps(x,\tau,0)\dd\tau
\]
which shows that $\gamma_\eps(x,t,0)$ is c-bounded and also concludes the proof that it belongs to $\cE_{M,L^0}(\Omega,\R^2)$. Therefore, if $v_{0\eps}$ belongs to
$\cE_{M,L^0}(\Omega,\R)$ then $v_\eps(x,t) = v_{0\eps}(\gamma_\eps(x,t,0))$ belongs to $\cE_{M,L^0}(\Omega,\R^2)$ and represents a pathwise solution to
\eqref{eq:transportG}. Uniqueness of the Colombeau solution can be shown by the same argument as in the proof of Proposition\;\ref{prop:tonly}. In conclusion, given
$v_0\in \cG_{L^0}(\Omega,\R)$, the stochastic transport equation \eqref{eq:transport} with coefficient \eqref{eq:transportcoeff} has a unique solution
$v\in \cG_{L^0}(\Omega,\R^2)$, as claimed.

Suppose that the process $X$ defining $c = \iota(X)$ has (almost surely) continuous paths which are of infinite length on every finite interval (as is the case with Brownian motion or the Ornstein-Uhlenbeck process,
\cite[Remark (5.2.8) and Section 8.3]{Arnold:1974}).
Then the relation \eqref{eq:arclength} entails that
\begin{equation}\label{eq:limchar}
   \lim_{\eps\to 0}\gamma_\eps(x,t,0) = x.
\end{equation}
Indeed, if \eqref{eq:limchar} does not hold, there is $x_0 < x$ and a sequence $\eps_k\to 0$ such that $\gamma_{\eps_k}(x,t,0) < x_0$ for all $k\in\N$.
Since $L_\eps(x)$ is a strictly increasing function of $x$, we have that $L_{\eps_k}(\gamma_{\eps_k}(x,t,0)) < L_{\eps_k}(x_0)$, and so \eqref{eq:arclength} implies that
$t > L_{\eps_k}(x) - L_{\eps_k}(x_0) = L_{\eps_k}$, the arclength of the curve $C_{\eps_k}$ between $x_0$ and $x$. We will obtain a contradiction by showing that this arclength tends to infinity as $\eps_k\to 0$. For $\omega$ in a set of full measure, we have that $X_{\eps_k} = X\ast(\chi\rho_{\eps_k}) \to X$ uniformly on
the interval $[x_0, x]$ (cf. Lemma\;\ref{lem:uniformconvergence}). Take a partition $x_0 < x_1 < \ldots < x_n = x$ of this interval. Given $\eta > 0$, we have that $|X_{\eps_k}(x_i) - X(x_i)| \leq \eta$ for all $k \geq k_0(\eta)$ and all $i$.
Let $L_{\eps_k,n}$ be the length of the polygonal chain through the points $(x_i, X_{\eps_k}(x_i))$. Then
\[
  L_{\eps_k,n} \geq \sum_{i=1}^n |X_{\eps_k}(x_i) - X_{\eps_k}(x_{i-1})| \geq \sum_{i=1}^n |X(x_i) - X(x_{i-1})| - 2n\eta = TV_n - 2n\eta,
\]
where $TV_n$ converges to the total variation of $X$ on $[x_0,x]$, which is infinite. Given $M > 0$, choose $n$ such that $TV_n \geq M$ and take $\eta = 1/n$ with corresponding $k_0(\eta)$. Then $L_{\eps_k} \geq L_{\eps_k,n} \geq M-2$ for
$k \geq k_0(\eta)$. Suitably choosing $M$, this contradicts that $L_{\eps_k} < t$. Thus \eqref{eq:limchar} holds.

If in addition $v_{0\eps} = v_0^\circ \ast (\chi\rho_\eps)$ for some continuous function $v_0^\circ$, then
\[
    \lim_{\eps\to 0} v_\eps(x,t) = v_0^\circ(x)
\]
almost surely for every $(x,t)$. Thus waves do not propagate on a curve of locally infinite length.

On the other hand, if $X$ has continuously differentiable paths, $X \in L^0(\Omega:\cC^1(\R))$, then $c_\eps'(x)$ converges locally uniformly to $X'(x)$, and so $L_\eps(z)$ converges locally uniformly to $L(z) = \int_0^z\sqrt{1+ X'(y)^2}\dd y$, which obviously is a strictly increasing, locally bi-Lipschitz function on $\R$. Using the strict monotonicity of $L$, the Lipschitz property of $L^{-1}$ and invoking \cite[Theorem 1]{Barvinek:1991} it follows
that $L_\eps^{-1}$ converges locally uniformly as well, resulting in the limiting characteristic curves
\[
   \lim_{\eps\to 0}\gamma_{\eps}(x,t,0) = \lim_{\eps\to 0} L_\eps^{-1}(L_{\eps}(x) - t) = L^{-1}(L(x) - t)) =: \gamma(x,t,0).
\]
Again, if in addition $v_{0\eps} = v_0^\circ \ast (\chi\rho_\eps)$ for some continuous function $v_0^\circ$, then
\[
    \lim_{\eps\to 0} v_\eps(x,t) = v_0^\circ(\gamma(x,t,0)),
\]
the usual propagation behavior.

{\bf Back to the wave equation.} Applying these results to the wave equation \eqref{eq:GeometricWave}, observe that the corresponding system \eqref{eq:geometricwavesyst} has the characteristic curves
\begin{equation}\label{eq:wavecharacteristics}
  \gamma_{\eps}^+(x,t,0) =  L_\eps^{-1}(L_{\eps}(x) - t),\qquad \gamma_{\eps}^-(x,t,0) =  L_\eps^{-1}(L_{\eps}(x) + t).
\end{equation}
From \eqref{eq:geometricwavesyst} we have that
\[
u^\pm_\eps(x,t) = u_{1\eps}(\gamma_\eps^\pm(x,t,0)) \mp \lambda_\eps(\gamma_\eps^\pm(x,t,0))u_{0\eps}'(\gamma_\eps^\pm(x,t,0))
\]
and $u_\eps(x,t) = u_{0,\eps}(x) + \frac12\int_0^t (u_\eps^+(x,s) + u_\eps^-(x,s))\dd s$.
Using the calculation \eqref{eq:dtgamma}, it follows from \eqref{eq:wavecharacteristics} that
\[
  \frac{\p}{\p t}\gamma_\eps^\pm(x,t,0) = \mp\lambda_\eps(\gamma_\eps^\pm(x,t,0)).
\]
Thus evaluating the integral and applying the substitution rule finally shows that $u_\eps(x,t)$ is given by
\[
   \frac12\Big(u_{0\eps}(\gamma_\eps^+(x,t,0)) + u_{0\eps}(\gamma_\eps^-(x,t,0))
      + \int_0^t \big(u_{1\eps}(\gamma_\eps^+(x,s,0)) + u_{1\eps}(\gamma_\eps^-(x,s,0))\big)\dd s\Big)
\]
from where the limiting behavior similar to the single transport equation can be read off.
\end{ex}

\begin{ex}\label{ex:Ogawa}
This example resumes Ogawa's analysis \cite{Ogawa} of the transport equation with white noise in time as coefficient. Thus we study the equation
\begin{equation}\label{eq:Ogawa}
(\p_t + \dot{W}(t)\p_x) u(x,t) = 0,\qquad u(x,0) = u_0(x)
\end{equation}
where $\dot{W}$ is Gaussian white noise in time. (For simplicity, we drop the additional deterministic terms in \cite[Equation (1.1)]{Ogawa} which would not substantially change the discussion.)
We take $u_0\in \cG_{L^0}(\Omega,\R)$ and interpret the coefficient as $\iota(\dot{W})\in \cG_{L^0}(\Omega,\R)$. Let $W$ be a one-parameter Brownian motion (Example\;\ref{ex:Brownian}(a)) and  $W_\eps = W\ast(\chi\rho_\eps)$ be a representative of $\iota(W)$. A representative
of $\iota(\dot{W})$ is given by $\dot{W}_\eps(t,\omega) = \frac{\dd}{\dd t}W_\eps(t,\omega)$. Then
\[
    \int_\tau^t \dot{W}_\eps(s,\omega)\dd s = W_\eps(t,\omega) - W_\eps(\tau,\omega);
\]
this shows that $\iota(\dot{W})$ is of almost certain local $L^1$-type. Thus Proposition\;\ref{prop:tonly} entails that \eqref{eq:Ogawa} has a unique solution $u\in \cG_{L^0}(\Omega,\R^2)$, represented by
\begin{equation}\label{eq:Ogawasolu}
   u_\eps(x,t,\omega) = u_{0\eps}\big(x - W_\eps(t,\omega),\omega\big).
\end{equation}
For a further analysis of its properties we focus on positive times.
First recall that $\rE(W_\eps(t)) = 0$. Indeed,
\begin{eqnarray*}
\rE(W_\eps(t)) &=& \int_\Omega W_\eps(t,\omega)\dd P(\omega) = \int_\Omega\int_{-\infty}^\infty\chi(t-s)\rho_\eps(t-s)W(s,\omega) \dd s \dd P(\omega)\\
&=& \int_{-\infty}^\infty\chi(t-s)\rho_\eps(t-s)\int_\Omega W(s,\omega) \dd P(\omega)\dd s\\
& =& \int_{-\infty}^\infty\chi(t-s)\rho_\eps(t-s)\rE(W(s))\dd s = 0.
\end{eqnarray*}
The transition from the first to the second line is valid due to Fubini's theorem, using that $W \in\cC(\R:L^1(\Omega))$ and the $s$-integration is over a compact interval.
A similar argument using the fact that $W \in\cC(\R:L^2(\Omega))$ shows that the variance of $W_\eps(t)$ can be computed as
\begin{equation}\label{eq:sigmaeps}
  \rE(W_\eps(t)^2) = \int_{-\infty}^\infty \int_{-\infty}^\infty \chi(t-s)\rho_\eps(t-s)\chi(t-s')\rho_\eps(t-s')\psi(s,s')\dd s\dd s',
\end{equation}
where $\psi(s,s') = \rE(W(s)W(s'))$, i.e, 
\[
\rE(W_\eps(t)^2) = \big((\chi\rho_\eps\otimes\chi\rho_\eps)\ast \psi\big)(t,t).
\]
Since $(\chi\rho_\eps\otimes\chi\rho_\eps)\ast \psi = ((\chi\otimes\chi)(\rho_\eps\otimes\rho_\eps))\ast \psi$ converges locally uniformly to the continuous function $\psi$
(Lemma\;\ref{lem:uniformconvergence}), it follows that for $t$ in a compact subinterval $I\subset (0,\infty)$, $\rE(W_\eps(t)^2)$ converges uniformly to $\psi(t,t) = t$ and hence is strictly bounded away from zero
as $\eps\to 0$.

Note that when $t$ varies in the compact interval $I$, the domain of integration in \eqref{eq:sigmaeps} is bounded; we may differentiate under the integral sign, and any $t$-derivative
of $\rE(W_\eps(t)^2)$ is bounded by a negative power of $\eps$ on $I$. The boundedness from below on $I$ shows that
\[
\sigma_\eps(t) = \sqrt{\rE(W_\eps(t)^2)}
\]
belongs to $\cE_M(0,\infty)$ and defines an invertible element of $\cG(0,\infty)$. Further, for $t>0$, $W_\eps(t)$ has the non-degenerate Gaussian distribution $\cN(0, \sigma_\eps(t)^2)$, whose density
\[
p_\eps(x,t) = \frac{1}{\sqrt{2\pi}\sigma_\eps(t)}\exp\Big(-\frac{x^2}{2\sigma_\eps(t)^2}\Big)
\]
also defines an element of $\cG(\R\times(0,\infty))$. In addition, for $0<t_0<t_1$ and $\alpha,\beta\in\N_0$,
\begin{equation}\label{eq:L1moderateness}
\sup_{t\in[t_0,t_1]}\|\p_x^\alpha\p_t^\beta p_\eps(\cdot,t)\|_{L^1(\R)} \in \cE_M.
\end{equation}
Our next task will be to study the expectation value of the solution $u\in \cG_{L^0}(\Omega,\R^2)$ to \eqref{eq:Ogawa}. Theorem\;\ref{thm:gensolup} is not applicable; nevertheless, the expectation value of $u$ exists as an element oF $\cG(\R\times(0,\infty))$ in certain special cases.

Let $u_0\in \cG(\R)$ be a uniformly bounded (deterministic) Colombeau function, i.e., for some representative $(u_{0\eps})_{\eps\in(0,1]}$, $\sup_{x\in\R}|u_{0\eps}(x)|$ is bounded independently of $\eps$. Then
\[
  \rE(u_\eps(x,t)) = \rE\big(u_{0\eps}(x-W_\eps(t)\big) = \int_{-\infty}^\infty u_{0\eps}(x-y)p_\eps(y,t)\dd y = \big(u_{0\eps}\ast p_\eps(\cdot,t)\big)(x).
\]
We have
\[
  \sup_{t\in[t_0,t_1]}\|\p_x^\alpha\p_t^\beta (u_{0\eps}\ast p_\eps(\cdot,t))\|_{L^\infty(\R)} \leq \|u_{0\eps}\|_{L^\infty(\R)}\sup_{t\in[t_0,t_1]}\|\p_x^\alpha\p_t^\beta p_\eps(\cdot,t)\|_{L^1(\R)}
\]
By \eqref{eq:L1moderateness}, the net $\big(\rE(u_\eps(x,t))\big)_{\eps\in(0,1]}$ belongs to $\cE_M(\R\times(0,\infty)$ and may serve to define the expectation value
\[
  \overline{u}(x,t) = \rE(u(x,t)) \in \cG(\R\times(0,\infty)).
\]
We now show that Ogawa's observation on the averaged quantity in \cite[Proposition 5]{Ogawa} has its counterpart in our setting. 

Let us start with classical (deterministic) initial data $u^\circ_0 \in \cD'(\R)$. 

\emph{The setting of Ogawa.} If $u^\circ_0$ is a continuously differentiable function with bounded derivative, then \cite[Theorem]{Ogawa} states that equation \eqref{eq:Ogawa} has a unique solution
in the sense of \cite[Section 3]{Ogawa} with probability one, and it is given by $u^\circ(x,t,\omega) = u^\circ_0(x-W(t,\omega))$. Further, if the function $\overline{u^\circ}(x,t) = \rE(u^\circ(x,t))$ 
is sufficiently regular, then it satisfies the heat equation (\cite[Proposition 5]{Ogawa})
\[
   \p_t \overline{u^\circ} - \frac12\p_x^2 \overline{u^\circ} = 0,\qquad \overline{u^\circ}(x,0) = u^\circ_0(x).
\]
\emph{The Colombeau setting.} Without further regularity assumptions on $u^\circ_0 \in \cD'(\R)$, we can take initial data $u_0 = \iota(u^\circ_0) \in \cG(\R)$, and equation \eqref{eq:Ogawa} will have a unique 
solution $u\in \cG_{L^0}(\Omega,\R^2)$. As was shown above, it is represented by $u_\eps(x,t,\omega) = u_{0\eps}\big(x - W_\eps(t,\omega)\big)$.
Further, if $u^\circ_0 \in L^\infty(\R)$, then $u_0 = \iota(u^\circ_0)$ is a uniformly bounded Colombeau function, the expectation value $\overline{u}(x,t) = \rE(u(x,t))$ exists as an element of $\cG(\R\times(0,\infty))$
and it satisfies the heat equation
\[
   \p_t \overline{u} - \frac12\p_x^2 \overline{u} = 0 \quad{\rm in}\quad \cG(\R\times (0,\infty)).
\]
(This follows immediately from the fact that the kernels $p_\eps(t,x)$ satisfy the heat equation.)

\emph{Transition between the two solution concepts.} Assume that $u^\circ_0$ is a uniformly continuous, bounded function on $\R$. A slight modification of Lemma\;\ref{lem:uniformconvergence} shows that 
$\|u_{0\eps} - u^\circ_0\|_{L^\infty(\R)} \to 0$ as $\eps\to 0$. Since $W_\eps(t,\omega)\to W(t,\omega)$ locally uniformly, we have that, for almost all $\omega$,
\[
   u_\eps(x,t,\omega) = u_{0\eps}\big(x - W_\eps(t,\omega)\big) \to u^\circ_0(x-W(t,\omega)) = u^\circ(x,t,\omega)
\]
locally uniformly on $\R\times [0,\infty)$, that is, the Colombeau solution $u$ is associated with Ogawa's solution $u^\circ$.

Concerning the expectation values we have, where $p(x,t)$ denotes the probability density of $\cN(0,t)$,
\begin{eqnarray*}
u_{0\eps}\ast p_\eps(\cdot,t) - u^\circ_0\ast p(\cdot,t) = (u_{0\eps} - u^\circ_0)\ast p_\eps(\cdot,t) + u^\circ_0\ast(p_\eps(\cdot,t) - p(\cdot,t)).
\end{eqnarray*}
Both summands converges to zero in $L^\infty(\R)$; the first one, because the $L^1$-norm of $p_\eps(\cdot,t)$ equals one, the second one, because
$p_\eps(\cdot,t) - p(\cdot,t)$ converges to zero in $L^1$. This shows that the expectation value of the Colombeau solution $\overline{u}$ is associated with
the expectation value $\overline{u^\circ}$ of Ogawa's solution.
\end{ex}

%
%
\appendix
\section{Appendix}
\label{Sec:appendix}
%
%
The appendix serves to (i) recall details on the embedding of distributions and generalized stochastic processes; (ii) provide some examples of such processes and
highlight how the required properties of the embedded versions can be guaranteed; (iii) substantiating the claims of Remark\;\ref{rem:inclusions} by means of counterexamples.

Properties of the embedding $\iota:\cD'(O) \to \cG(O)$ on open subsets $O\subset\R^d$ have been elaborated in detail in \cite[Section 1.2.2]{GKOS}. We wish to
show here that the simple form of \eqref{eq:imbedding} leads to the same embedding on $O=\R^d$. Also, in view of the stochastic generalization, it appears useful
to explicitly state some of the relevant estimates.

\emph{Embedding of distributions.} The embedding of the space of Schwartz distributions in $\cG(\R^d)$ has been defined by equation \eqref{eq:imbedding}
with a fixed mollifier of the form \eqref{eq:molli}. Here is the argument why it is well-defined. First, one has to show that the representative $\big(v \ast (\chi\rho_\eps)\big)_{\eps \in (0,1]}$ belongs to $\cE_M(\R^d)$. Let $v\in\cD'(\R^d)$, $K$ a compact subset of $\R^d$ and $\alpha$ a multi-index.
For $x\in K$, the support of the test function $y\to \chi(x-y)\rho_\eps(x-y)$ is contained in the compact set $L = \{x\} - \supp \chi$. By the continuity of the functional $v$ on $\cD(\R^d)$, there are $C>0$ and $m\in\N$ such that
\[
  |\langle \p^\alpha v(y), \chi(x-y)\rho_\eps(x-y)\rangle| \leq C \sup_{y\in L}\sup_{|\beta|\leq m}|\p^{\alpha+\beta}(\chi(x-y)\rho_\eps(x-y))|.
\]
Thus
\begin{equation}\label{eq:continuity}
   \sup_{x\in K}\big|\p^\alpha v\ast(\chi\rho_\eps)(x)\big| \leq C\sup_{z\in K-L}\sup_{|\beta|\leq m}|\p^{\alpha+\beta}(\chi(z)\rho_\eps(z))|.
\end{equation}
The latter term is obviously bounded by a negative power of $\eps$, establishing moderateness of $\big(v \ast (\chi\rho_\eps)\big)_{\eps \in (0,1]}$.

To prove the injectivity of $\iota$, take $\psi\in\cD(\R^d)$ and note first that $(\check{\chi}\check{\rho}_\eps)\ast\psi \to \psi$ in $\cD(\R^d)$, where the check denotes inflection. Therefore,
\begin{equation}\label{eq:injectivity}
  \langle v\ast (\chi\rho_\eps),\psi\rangle = \langle v, (\check{\chi}\check{\rho}_\eps)\ast\psi\rangle \to \langle v, \psi \rangle.
\end{equation}
On the other hand, if $K = \supp \psi$ and $(v\ast (\chi\rho_\eps))_{\eps \in (0,1]} \in \cN(\R^d)$, then $\sup_{x\in K}|v\ast(\chi\rho_\eps)(x)|\leq C_q\eps^q$ for whatever $q > 0$, so
$\langle v \ast (\chi\rho_\eps),\psi\rangle \to 0$. It follows that $\langle v,\psi\rangle = 0$.

We show next that the embedding coincides with the (usual) embedding as defined in \cite[Section 1.2.2]{GKOS}.
For this, it suffices to show that the two embeddings coincide on $\cE'(\R^d)$ (using \cite[Theorem 1.2.20]{GKOS}), that is,
\begin{equation}\label{eq:nocutoff}
     \iota(v) = \ \mbox{class of}\ \big(v \ast \rho_\eps\big)_{\eps \in (0,1]}
\end{equation}
for $v\in\cE'(\R^d)$. This in turn means that the net $(v \ast (\chi\rho_\eps - \rho_\eps))_{\eps \in (0,1]}$ belongs to $\cN(\R^d)$. Fixing $x\in\R^d$ and using that $v$ is a continuous linear form on $\cE(\R^d)$, there is $L\Subset \R^d$, $C>0$ and $m\in\N$
such that
\[
  |\langle v(y), (\chi(x-y)-1)\rho_\eps(x-y)\rangle| \leq C \sup_{y\in L}\sup_{|\alpha|\leq m}|\p^\alpha((\chi(x-y)-1)\rho_\eps(x-y))|.
\]
Thus, given $K\Subset\R^d$,
\[
  \sup_{x\in K}|\langle v(y), (\chi(x-y)-1)\rho_\eps(x-y)\rangle| \leq C \sup_{z\in K-L}\sup_{|\alpha\leq m}|\p^\alpha|((\chi(z)-1)\rho_\eps(z))|.
\]
The latter term is easily seen to be negligible, using that $\rho$ belongs to $\cS(\R^d)$ and that $\p^\beta(\chi(z) - 1) \equiv 0$ for $|z| \leq\eta$ for some $\eta > 0$ and all $\beta$.

It follows that the embedding \eqref{eq:imbedding} has all the properties elaborated in \cite[Section 1.2.2]{GKOS}. In particular, it does not depend on the choice of the cut-off $\chi$; it commutes with partial derivatives; it coincides with the embedding $\sigma$ from \eqref{eq:standardimbedding} on $\cC^\infty(\R^d)$. Finally, the cut-off is not needed when $v\in\cS'(\R^d)$, that is, $\iota$ is given by \eqref{eq:nocutoff} when $v\in \cS'(\R^d)$, \cite[Proposition 1.2.21]{GKOS}.

\emph{Embedding of generalized stochastic processes.} Embeddings of generalized stochastic processes ($\cD'(\R^d)$-valued processes and distributions with values in $L^p(\Omega)$) have been described in \cite{GOPS2018,Mirkov:2009}.
Using results of \cite{MO:2024}, a unified and somewhat simpler treatment is possible, which we adopt here.

As noted in Section\;\ref{Sec:Colombeaurandfun}, a random Schwartz distribution is a weakly measurable map $X: \Omega\to \cD'(\R^d)$, i.e., for all $\psi\in\cD(\R^d)$, the maps
$\omega\to \langle X(\omega),\psi\rangle$ are measurable \cite[Section 1.3(ii)]{Hida}, \cite{Ullrich:1957}. The space of random Schwartz distributions may be denoted by $L^0(\Omega:\cD'(\R^d))$ where $\cD'(\R^d)$ is equipped with the Borel $\sigma$-algebra generated by the weak topology. A \emph{linear random functional} or \emph{generalized random function} on $\cD(\R^d)$ is a sequentially continuous linear map $\cD(\R^d)\to L^0(\Omega)$, the latter space equipped with the metric topology of convergence in probability \cite[Section III.1.2]{GelfandVilenkin:1964}, \cite[Section 2.3]{Ito:1984}; we use the notation $\cD'(\R^d:L^0(\Omega))$. It is known that every linear random functional has a version which is a random Schwarz distribution \cite[Theorem 2.3.3]{Ito:1984}. However, more is true. The linear map $\kappa: L^0(\Omega:\cD'(\R^d)) \to \cD'(\R^d:L^0(\Omega))$, defined by $\kappa(X)(\varphi)(\omega) = \langle X(\omega), \varphi\rangle$ for $\varphi\in \cD(\R^d)$ and almost all $\omega$, is a bijection, as has been shown in \cite{MO:2024}. Denoting by $\cD'(\R^d:L^p(\Omega))$ the space of linear continuous maps $\cD(\R^d)\to L^p(\Omega)$, $p\geq 1$, we have the inclusions
\[
   \cD'(\R^d:L^p(\Omega)) \subset \cD'(\R^d:L^0(\Omega)) = L^0(\Omega:\cD'(\R^d)),
\]
and we can embed these spaces simultaneously into $\cG_{L^0}(\Omega,\R^d)$ and into $\cG_{L^p}(\Omega,\R^d)$, respectively.

We begin by embedding $L^0(\Omega:\cD'(\R^d))$. If $X$ is a random Schwartz distribution and $\omega\in\Omega$, we can form the net
$\big(X(\omega) \ast (\chi\rho_\eps)\big)_{\eps \in (0,1]}$. At fixed $\eps$ it is a random field on $\R^d$ given by
\[
   (x,\omega)\to \big(X(\omega)\ast(\chi\rho_\eps)\big)(x) = \langle X(y,\omega),\chi(x-y)\rho_\eps(x-y)\rangle
\]
and it obviously satisfies the conditions (S) and (L$_0$) of Section\;\ref{Sec:Colombeaurandfun}.
The estimate \eqref{eq:continuity} with $X(\omega)$ in place of $v$ shows that the net $\big(X(\omega) \ast (\chi\rho_\eps)\big)_{\eps \in (0,1]}$ is almost surely moderate, i.e., satisfies condition (M$_0$).
Going over to the factor spaces, we thus obtain an embedding
\begin{equation} \label{eq:randomimbedding}
   \iota:L^0(\Omega:\cD'(\R^d)) \to \cG_{L^0}(\Omega,\R^d);
\end{equation}
the injectivity follows from the same argument as in the deterministic case. This embedding shares all the properties of the embedding of $\cD'(\R^d)$ into $\cG(\R^d)$, pathwise at fixed $\omega\in\Omega$. In particular, it commutes with derivatives and coincides with the embedding
$\sigma: F\to {\rm class\ of\ } [(\eps,x,\omega) \to F(x,\omega)]$ for random fields $F$ with smooth paths.

Let $1\leq p \leq \infty$ and assume in addition that $\langle X,\varphi\rangle \in L^p(\Omega)$ for all $\varphi\in \cD(\R^d)$ and that the map $\varphi \to \langle X,\varphi\rangle$ is linear and continuous from $\cD(\R^d)$ to $L^p(\Omega)$, that is, $X\in \cD'(\R^d:L^p(\Omega))$. 
Thus at fixed $x\in\R^d$, $\alpha\in\N_0^d$ and $0<\eps\leq 1$, the map $\omega\to\big(X(\omega)\ast(\chi\rho_\eps)\big)(x)$ belongs to $L^p(\Omega)$; 
the continuity entails an estimate analogous to \eqref{eq:continuity}, namely
\[
  \sup_{x\in K}\big\|\p^\alpha X\ast(\chi\rho_\eps)(x)\big\|_{L^p(\Omega)} \leq C\sup_{z\in K-L}\sup_{|\beta|\leq m}|\p^{\alpha+\beta}(\chi(z)\rho_\eps(z))|,
\]
showing that conditions (L$_p$) and (M$_p$) hold. Finally, if $\big(X\ast(\chi\rho_\eps)(x)\big)_{\eps\in (0,1]}$ satisfies condition (N$_p$), then $X\equiv 0$, as shown by the same argument as in the deterministic case, observing that convergence in \eqref{eq:injectivity} takes place in $L^p(\Omega)$. This way we obtain the embedding
$\iota:\cD'(\R^d:L^p(\Omega)) \to \cG_{L^p}(\Omega,\R^d)$.

We next collect some examples of classical random fields and random Schwartz distributions which are relevant to the examples of Section\;\ref{Sec:apps}.

\begin{ex}\label{ex:continuousprocess}
(a) Every (classical) process $X:\R^d\times\Omega \to \R$ with (almost surely) continuous paths defines a random Schwartz distribution via
\begin{equation}\label{eq:continuousprocess}
\langle X(\omega),\varphi\rangle = \int_{\R^d}\varphi(x)X(x,\omega)\dd x.
\end{equation}
Indeed, $\langle X(\omega),\varphi\rangle$ is measurable as a limit of Riemann sums. If $\varphi_n\to \varphi$ in $\cD(\R^d)$, then clearly $\int\varphi_n(x)X(x,\omega)\dd x \to \int\varphi(x)X(x,\omega)\dd x$ almost surely.

(b) If $1\leq p \leq\infty$ and, in addition, $X$ belongs to $\cC(\R^d:L^p(\Omega))$, then \eqref{eq:continuousprocess} results in a process belonging to $\cD'(\R^d:L^p(\Omega))$. Indeed, processes with continuous paths are jointly measurable in $(x,\omega)$. We may apply Minkowski's inequality for integrals to obtain
\[
\left(\int_\Omega \big\vert\int_{\R^d}\varphi(x) X(x,\omega)\dd x\big\vert^p\dd P(\omega)\right)^{1/p}
   \leq \int_{\R^d}\left(\int_\Omega|X(x,\omega|^p\dd P(\omega)\right)^{1/p}|\varphi(x)|\dd x
\]
with a suitable modification for $p=\infty$. It follows that $\varphi\to \langle X,\varphi\rangle$ is a continuous map on $\cD(\R^d)$ with values in $L^p(\Omega)$.

(c) If $Y$ is the $\alpha$th distributional derivative of a process $X$ as in (a), (b), that is
\[
   \langle Y(\omega),\varphi\rangle  = (-1)^{|\alpha|} \langle X(\omega), \p^{\alpha}\varphi\rangle
\]
it clearly defines a process in $L^0(\Omega:\cD'(\R^d))$ as well as in $\cD'(\R^d:L^p(\Omega))$. The embedded process $\iota(Y)$ in $\cG_{L^p}(\Omega,\R^d)$ ($p=0$ or $p\geq 1$, respectively) is represented by the net
\[
   Y_\eps(x,\omega) = \int_{\R^d} \p^{\alpha}(\chi\rho_\eps)(x-y)X(y,\omega)\dd x.
\]
\end{ex}

\begin{ex}\label{ex:Brownian}
(a) Recall that a \emph{$d$-parameter Brownian motion} is a Gaussian stochastic process $W$ on $\R^d$ such that
(i) $W(0) = 0$ almost surely; (ii) $\rE(W(x)) = 0$ for all $x\in\R^d$; (iii) $E(W(x) W(y)) = \prod_{i=1}^d x_i\wedge y_i$ for $(x,y)\in \R_+^d\times \R_+^d$ with an appropriate sign-depending extension to all of $\R^d\times \R^d$
(see \cite[Section 2.1.1]{HoldOksUbZhang}); (iv) $W$ has almost surely continuous paths.

The special form of the covariance shows that $\rE\big((X(x+h) - X(x))^2\big) = \cO(|h|)$ as $h\to 0$. Thus $W$ belongs to $\cC(\R^d:L^2(\Omega))$, so that $W\in \cD'(\R^d:L^2(\Omega))\subset L^0(\Omega:\cD'(\R^d))$.

(b) Taking the derivative of order $d$ in the sense of distributions leads to
\emph{Gaussian white noise:}
\begin{equation}\label{eq:GaussianWhiteNoise}
\dot{W} = \p_{x_1}\ldots \p_{x_d} W,\qquad \langle \dot{W},\varphi\rangle = (-1)^d\langle W, \p_{x_1}\ldots \p_{x_d}\varphi\rangle
\end{equation}
Gaussian white noise is a Gaussian random Schwartz distribution, i.e., all the random variables $\langle\dot{W},\varphi_1\rangle, \ldots, \langle\dot{W},\varphi_n\rangle$ are jointly Gaussian for every choice of $\varphi_i\in \cD(\R^d)$, $n\in\N$.
It has mean zero and covariance functional
\begin{equation}\label{eq:ItoIsometry}
\rE(\langle\dot{W},\varphi\rangle\langle\dot{W},\psi\rangle) = \int_{\R^d}\,\varphi(y)\psi(y)\,dy,
\end{equation}
as follow from the fact that $\p_{x_1}\ldots \p_{x_d}\prod_{i=1}^d x_i\wedge y_i = \prod_{i=1}^d\delta(x_i - y_i)$ in the sense of distributions.
The relation \eqref{eq:ItoIsometry} constitutes the so-called It\^{o} isometry.
Clearly, $\dot{W}$ belongs to $\cD'(\R^d:L^2(\Omega))$.
\end{ex}

\begin{ex}\label{ex:randomfield}
A \emph{second order random field} on $\R^d$ is a family of random variables $\{X(x):x\in \R^d\}$ such that its variances $\rV(X(x))$ are finite for all $x\in\R^d$. Its
covariance function $C(x,y) = {\rm COV}(X(x),X(y))$ is well defined and nonnegative definite. Given a nonnegative definite function $C:\R^d\times\R^d$, there is always a probability space $\Omega$ and a Gaussian random field $X:\R^d\times \Omega\to \R$ which has $C(x,y)$ as its covariance function \cite[Section 37.1]{LoeveII:1977}.
Mean square continuity or differentiability can be read off from corresponding properties of the covariance function along the diagonal \cite[Section 37.2]{LoeveII:1977}.
If the covariance function $C(x,y)$ depends only on the distance $|x-y|$, the random field is \emph{isotropic}.
Widely used covariance functions of isotropic Gaussian random fields are
\[
   C(x,y) = \sigma^2\exp(|x-y|/\ell)\quad {\rm or}\quad \sigma^2\exp(|x-y|^2/\ell^2),
\]
where the parameter $\ell$ signifies correlation length.
In the first case, the paths are (almost surely) continuous \cite[Theorem 1.4.1]{Adler:2007}, but nowhere differentiable \cite[Theorem 6]{Cambanis:1973}. In the second case, the paths are smooth according to, e.g., \cite[Theorem B]{Kawata:1970}
or \cite[Theorem 1]{Scheurer:2010}.
\end{ex}

\begin{rem}\label{rem:boundedness}
This remark serves to show how the boundedness conditions of Definition\;\ref{def:bounded} can be enforced for Colombeau random fields, obtained as embeddedings
of classical or generalized stochastic processes.

(a) Enforcing boundedness: If $X$ is a Gaussian random field with mean zero and unit variance as in Example\;\ref{ex:randomfield}, one may transform it to a random field
\[
  Y(x) = F^{-1}(\Phi(X(x)))
\]
where $F$ is a given distribution function and $\Phi$ the distribution function of a standard Gaussian variable. Then the marginals $Y(x)$ are distributed according to $F$. If $F^{-1}$ has bounded range, then $Y(x)$ will be globally bounded. Further, $Y$ inherits the pathwise properties of $X$, provided $F$ is sufficiently smooth. The embedding \eqref{eq:nocutoff} produces an element of $\cG_{L^0}(\Omega:\R^d)$, which is globally uniformly bounded (in $x,\omega,\eps$). In particular, in the case $d=2$ it will be uniformly bounded on strips. For a detailed study of the properties of these \emph{translation processes} see \cite{Grigoriu:2007}.

Alternatively, if $X\in L^0(\Omega:\cD'(\R^d))$ one may first embed $X$ into $\cG_{L^0}(\Omega:\R^d)$ via \eqref{eq:randomimbedding} and then multiply $\iota (X)$ by a smooth cut-off function
\[
\psi(y) = \left\{ \begin{array}{ll}
          M, & y \geq M\\
          y, & m+\eta < y < M-\eta\\
          m, & y \leq m
          \end{array}\right.
\]
to enforce the range $[m,M]$, where $m < M$ and $\eta > 0$ is small. This results in an element of $\cG_{L^0}(\Omega:\R^d)$ which again is globally uniformly bounded.

(b) Enforcing almost certain local logarithmic type: Replacing the mollifier \eqref{eq:molli} in the embedding \eqref{eq:imbedding}, \eqref{eq:randomimbedding} by a suitably scaled
mollifier $\rho_{\eta(\eps)}$ produces asymptotic bounds in terms of negative powers of $\eta(\eps)$. In general, if $X\in L^0(\Omega:\cD'(\R^d))$, the choice of
$\eta(\eps) = \log|\log\eps|$ will give that $\iota(X)\in \cG_{L^0}(\Omega:\R^d)$ is of almost certain local logarithmic type, together with all its derivatives.
Depending on the local order of $X(\omega)$, less stringent choices are possible. For example, if $X\in L^0(\Omega:\cD'(\R^d))$ has almost surely continuous paths, then
$\eta(\eps) = |\log\eps|$ suffices for $\iota(X)$ to be of almost certain local logarithmic type.
\end{rem}

\begin{ex}\label{ex:counterexamples}
Here are two counterexamples illustrating the non-embedding properties alluded to in Remark\;\ref{rem:inclusions}.

(a) The map $\cG_{L^{p'}}(\Omega,\R^d) \to \cG_{L^{p}}(\Omega,\R^d)$ need not be injective. This happens already in the case of Colombeau random variables.
We display an element of $\cE_{M,L^{p'}}(\Omega)$ which belongs to $\cN_{L^{p}}(\Omega)$ for all $p<p'$, but not to $\cN_{L^{p'}}(\Omega)$.
We take $\Omega = [0,\infty)$ with the standard exponential distribution. We set
\[
   u_\eps(\omega) = \left\{\begin{array}{ll}
       0, & \omega < \ee^{p'/\eps},\\[2pt]
       \ee^{1/\eps}, & \omega \geq \ee^{p'/\eps}.
       \end{array}\right.
\]
In particular, at fixed $\omega\geq \ee^{p'/\eps}$, $u_\eps(\omega) = 0$ for $0 < \eps \leq \eps(\omega) = p'/\log\omega$. (In particular, this shows that
the net $(u_\eps)$ defines an element of $\cN_{L^{0}}(\Omega)$.) Let $p\geq 1$. Then
\[
 \rE(u_\eps^p) = \ee^{p/\eps}\int_{\ee^{p'/\eps}}^\infty \ee^{-\omega}\dd\omega = \ee^{(p-p')/\eps},
\]
which is negligible for $p < p'$ but not negligible (namely equal to 1) for $p=p'$.

(b) $\cE_{M,L^1}(\Omega) \not\subset \cE_{M,L^0}(\Omega)$.
Indeed, take $\Omega = [-1,1]$ with the uniform distribution and set
\[
   v_\eps(\omega) = \ee^{1/\eps}\eins_{(-\ee^{-1/\eps},\ee^{-1/\eps})}(\omega), \qquad u_\eps(\omega) = v_\eps\big(\omega - \sin\tfrac1{\eps}\big).
\]
Then
\[
\rE(|u_\eps|) = \tfrac12\int_{-1}^1 u_\eps(\omega)\dd\omega \leq \tfrac12\int_{-1}^1 v_\eps(\omega)\dd\omega = 1,
\]
so $(u_\eps)_{\eps\in I} \in \cE_{M,L^1}(\Omega)$. On the other hand, for each fixed $\omega\in[-1,1]$, there is a sequence $\eps_k\to 0$ such that
$\omega = \sin\frac{1}{\eps_k}$ for all $k\in\N$. Then
\[
   u_{\eps_k}(\omega) = v_{\eps_k}\big(\omega - \sin\tfrac{1}{\eps_k}\big) = v_{\eps_k}(0) = \ee^{1/\eps_k},
\]
so $(u_\eps)_{\eps\in(0,1]} \not\in \cE_{M,L^0}(\Omega)$.
\end{ex}

\end{document}